\documentclass{amsart}
\usepackage{amssymb,latexsym}
\usepackage{amscd,amsthm}

\newtheorem{theorem}{Theorem}[section]
\newtheorem{lemma}[theorem]{Lemma}
\newtheorem{proposition}[theorem]{Proposition}
\newtheorem{corollary}[theorem]{Corollary}

\theoremstyle{definition}
\newtheorem{definition}[theorem]{Definition}
\newtheorem{fact}[theorem]{Fact}
\newtheorem*{remark}{Remark}

\newtheorem*{example}{Example}

\DeclareMathOperator{\Ext}{Ext} \DeclareMathOperator{\Hom}{Hom}
\DeclareMathOperator{\Tor}{Tor}
\DeclareMathOperator{\colim}{colim}
 \DeclareMathOperator{\cok}{cok}

\newcommand{\RHom}{\Hom_{R}}                 

\newcommand{\tensor}{\otimes}

\newcommand{\cat}[1]{\mathcal{#1}}           

\newcommand{\class}[1]{\mathcal{#1}}   

\newcommand{\Z}{\mathbb{Z}}

\newcommand{\mathcolon}{\colon\,} 
\newcommand{\ar}{\xrightarrow{}} 

\newcommand{\ch}{\textnormal{Ch}(R)}

\newcommand{\qcox}{\textnormal{Qco}(X)}
\newcommand{\chqcox}{\textnormal{Ch}(\textnormal{Qco}(X))}

\newcommand{\lclass}[1]{\widetilde{\class{#1}}}
\newcommand{\rclass}[1]{\widetilde{\class{#1}}}
\newcommand{\ldgclass}[1]{dg\widetilde{\class{#1}}}
\newcommand{\rdgclass}[1]{dg\widetilde{\class{#1}}}

\newcommand{\rightperp}[1]{#1^{\perp}}
\newcommand{\leftperp}[1]{{}^\perp #1}

\begin{document}

\title{A Quillen approach to Derived categories and tensor products}

\date{\today}

\author{James Gillespie}
\address{4000 University Drive\\
         Penn State McKeesport\\
         McKeesport, PA 15132-7698}
\email[Jim Gillespie]{jrg21@psu.edu}

\begin{abstract}
We put a monoidal model category structure (in the sense of
Quillen) on the category of chain complexes of quasi-coherent
sheaves over a quasi-compact and semi-separated scheme $X$. The
approach generalizes and simplifies the method used by the author
in~\cite{gillespie} and~\cite{gillespie-sheaves} to build monoidal
model structures on the category of chain complexes of modules
over a ring and chain complexes of sheaves over a ringed space.
Indeed, much of the paper is dedicated to showing that in any
Grothendieck category $\cat{G}$, a nice enough class of objects
$\class{F}$ (which we call a Kaplansky class) induces a model
structure on the category Ch($\cat{G}$) of chain complexes. We
also find simple conditions to put on $\class{F}$ which will
guarantee that our model structure in monoidal. We see that the
common model structures used in practice are all induced by such
Kaplansky classes.
\end{abstract}

\maketitle

\section{Introduction}\label{sec-introduction}

The classical theory of derived functors was introduced in Cartan
and Eilenberg's \emph{Homological Algebra}. The idea is now very
basic and is of fundamental importance to many branches of
mathematics. As we explain below, the author feels that the best
foundation for defining and treating derived functors is to use
Quillen's idea of a model structure.

The definition of derived functor appearing in Cartan and
Eilenberg's book depended on the existence of projective and
injective resolutions. Often however, especially in algebraic
geometry, one is usually dealing with a category in which
projective resolutions do not exist. For example, it is well-known
that we do not have projective resolutions in the various types of
sheaf categories. Therefore, Grothendieck and his school
generalized the definition of derived functor as to make no
mention of projective or injective objects at all. This definition
depends on first defining the derived category of an abelian
category $\cat{A}$. This is the ``category'' $\cat{D}(\cat{A})$
obtained by first forming the associated category Ch($\cat{A}$) of
\emph{unbounded} chain complexes, and then formally inverting the
homology isomorphisms. (The word ``category'' is in quotes because
it is only known in particular cases that the class of maps
between two objects in $\cat{D}(\cat{A})$ actually form a set.)
The derived category of an abelian category $\cat{A}$ is now
widely accepted as the proper setting to define and study derived
functors of all sorts. However, in general $\cat{D}(\cat{A})$ is
not at all easy to understand. Even the morphism sets seem
mysterious, without further analysis in each case.

However, Quillen's notion of a model category, appearing
in~\cite{quillen}, gave us a language and theory designed to deal
with categories exactly like derived categories. That is,
categories which are obtained by localizing with respect to a
class of morphisms. Indeed, with a model structure on
$\textnormal{Ch}(\cat{A})$, the morphism set
$\cat{D}(\cat{A})(X,Y)$ for complexes $X,Y \in \cat{D}(\cat{A})$
will be isomorphic to
$$\textnormal{Ch}(\cat{A})(QX,RY)/\sim$$ where $QX$ and $RY$ are
some sort of resolutions of $X$ and $Y$ respectively (called
``cofibrant'' and ``fibrant'' replacements, respectively) and
$\sim$ is the relation of chain homotopic maps. This alone gives
us a better understanding of the derived category, not to mention
the versatility that different model structures correspond to
different cofibrant and fibrant replacements. As an illustration,
take $\cat{A}$ to be the category of $R$-modules where $R$ is a
commutative ring with identity. There is a ``projective'' model
structure on $\ch$ in which cofibrant replacements correspond to
projective resolutions. There is also an ``injective'' model
structure where the fibrant replacements correspond to injective
coresolutions. Both these model structures are described in detail
in~\cite{hovey-model-categories}. The important result of
Spaltenstein in~\cite{spaltenstein} which says that every
unbounded complex has an injective resolution is automatic with
the existence of the injective model structure.

Quillen's theory also includes an easy way to verify the existence
of derived functors and also to \emph{compute} derived functors
using the cofibrant and fibrant approximations. If we look at the
$\ch$ example again, both the projective model structure and the
injective model structure are suitable to prove the existence and
compute the derived functors $\Ext^n_R(A,B)$. The first
corresponds to the classical computation of $\Ext^n_R(A,B)$ by
taking a projective resolution of $A$ and the second corresponds
to the classical definition of computing $\Ext^n_R(A,B)$ by taking
a coresolution of $B$.

However, not all model structures are suitable for studying all
derived functors. For example, the projective model structure on
$\ch$ is suitable to study $\Tor^R_n(A,B)$ but the injective model
structure is not. In general, if the category $\cat{A}$ has a
tensor product then in model category language we say that the
model structure on $\textnormal{Ch}(\cat{A})$ is \emph{monoidal}
if it is compatible with the tensor product on
$\textnormal{Ch}(\cat{A})$. It is now abundantly clear: If
$\cat{G}$ is a category without projective resolutions, what model
structure allows for an easy treatment of the derived tensor
product? In particular, we would like an answer for this when
$\cat{G}$ is the category $\qcox$ of quasi-coherent
$\class{O}_X$-modules on a scheme $(X, \class{O}_X)$. Clearly, the
literature suggests that there should exist a ``flat model
structure'' on Ch($\qcox$), which will allow computing
$\Tor^{\cat{G}}_n(A,B)$ via some sort of ``flat resolution''.

With this goal in mind, the author first constructed an analogous
model structure on $\textnormal{Ch}(\cat{G})$ when $\cat{G}$ is
the category of modules over a commutative ring $R$
in~\cite{gillespie}. This was generalized to the case of when
$\cat{G}$ is the category of sheaves of $\class{O}_X$-modules
where $(X, \class{O}_X)$ is any ringed space
in~\cite{gillespie-sheaves}. Through this experience and while
working on the problem when $\cat{G} = \qcox$ the author developed
a powerful theorem for building such homological model structures.
This method appears as Theorem~\ref{them-model structure}. The
theorem basically says that a suitable class of objects
$\class{F}$ in a Grothendieck category $\cat{G}$ gives rise to a
model structure on Ch($\cat{G}$). The theorem is very general and
in fact we show in Section~\ref{sec-other applications} that every
(homological) model structure the author has dealt with arises in
this way.

In Theorem~\ref{them-monoidal condition} we provide conditions
which will guarantee that the class $\class{F}$ induces a monoidal
model structure. In Section~\ref{sec-quasi-coherent sheaves} we
see that the class of flat objects in $\qcox$ induces a monoidal
model structure on Ch($\qcox$) when $X$ is a quasi-compact and
semi-separated scheme. Results from~\cite{lipman}
and~\cite{enochs-quasi-coherent} are essential and are used in the
construction of the flat model structure.

We point out that with all of the current literature on model
categories our flat model structure on $\textnormal{Ch}(\qcox)$
provides a strong foundation for doing homological algebra in
$\qcox$. For example, from~\cite{may3} it is now automatic that
$\cat{D}(\qcox)$ is triangulated in a way that is strongly
compatible with the derived tensor product. We also automatically
have the additivity property of generalized trace maps as defined
in~\cite{may3}.

We now summarize the layout of this paper. Our goal is to show
that a nice enough Kaplansky class $\class{F}$ in a Grothendieck
category $\cat{G}$ induces a (monoidal) model structure. We do
this by using Theorem~2.2 of~\cite{hovey} relating cotorsion pairs
and model structures. In Section~\ref{sec-small cotorsion pairs}
we develop some theory on small cotorsion pairs. The word
``small'' was first used to describe a cotorsion pair
in~\cite{hovey}. An analogy is that \emph{small} is to cotorsion
pair as \emph{cofibrantly generated} is to model category. The
material in this section was surely known to Hovey, but to the
author's knowledge, has not been written down in as much detail.
The material in Subsection~\ref{subsec-an associated cotorsion
pair of chain complexes is small} is new and is due to the author.
In particular, Proposition~\ref{prop-induced cotorsion pair of
complexes is small} greatly simplifies the problem of building the
flat model structure and makes great chunks of the work
in~\cite{gillespie} and~\cite{gillespie-sheaves} seem obsolete. In
Section~\ref{sec-locally cogenerated and kaplansky} we introduce
the notion of Kaplansky class in a Grothendieck category
$\cat{G}$. Our definition is a categorical version of the
definition given by Edgar Enochs
in~\cite{enochs-kaplansky-classes}. Enochs has proved the
existence of flat covers and cotorsion envelopes in several
algebraic categories, essentially by showing that the flat objects
form a Kaplansky class. We see in Section~\ref{sec-locally
cogenerated and kaplansky} that Kaplansky classes give rise to
cotorsion pairs and the cotorsion pairs give rise to homological
model structures. In Section~\ref{sec-monoidal-model structures}
we show that if the Kaplansky class satisfies certain
compatibilities with the tensor product, then the induced model
structure is monoidal in the sense
of~\cite{hovey-model-categories}. In
Section~\ref{sec-quasi-coherent sheaves} we prove the existence of
the flat model structure on Ch($\qcox$) when $X$ is a
quasi-compact, semi-separated scheme and in Section~\ref{sec-other
applications} we show that all of the usual homological model
structures come from Kaplansky classes.

The author would like to thank Mark Hovey for always answering
questions related to this work. A couple of lemmas which are
entirely due to Hovey are pointed out in the text. Thankyou also
to Edgar Enochs, and to Sergio Estrada and Leo Alonso Tarr\'io for
helping me understand a few things they already knew about
quasi-coherent sheaves.


\section{Preliminaries}\label{sec-preliminaries}

Prerequisites for this paper are a basic understanding of
Grothendieck categories, chain complexes, model categories,
cotorsion pairs, and (quasi-coherent) sheaves. The basic reference
used here for Grothendieck categories is~\cite{stenstrom}. The
author has used the texts~\cite{hartshorne} and (to a lesser
degree)~\cite{litaka} as references for sheaves, schemes and
quasi-coherent sheaves. Also see Appendix~B
of~\cite{thomason-trobaugh} for more advanced topics such as the
definition and basic facts concerning semi-separated schemes. The
author uses~\cite{hovey-model-categories} for referencing facts on
model categories.

Our work in this paper heavily rests on the work in~\cite{hovey}
and~\cite{gillespie}. The first paper laid out the basic interplay
between cotorsion pairs and model categories, while the second
focused exclusively on the interplay between cotorsion pairs of
chain complexes and (homological) model structures. In particular
we will use definitions and basic results from Section~3
of~\cite{gillespie}. There is a warning: Definition~3.11
in~\cite{gillespie} only makes sense if our category has enough
projectives and injectives. This mistake actually leaves a gap in
the construction of the flat model structure on
Ch($\class{O}_X$-Mod) which appeared in the
sequel~\cite{gillespie-sheaves}. The mistake is easily fixed.
Indeed one should replace the ``hereditary'' hypothesis with the
assumption appearing as condition (4) in Theorem~\ref{them-model
structure} of this paper. In this way we have avoided using the
word ``hereditary'' at all in a category without enough
projectives. In any case, we show in Section~\ref{sec-other
applications} that the flat model structure of
\cite{gillespie-sheaves} is just another corollary of our
Theorem~\ref{them-model structure}.

We also use language associated with locally $\kappa$-presentable
categories throughout much of the paper. This language along with
all the results we will need are summarized in
Appendix~\ref{appendix-locally-presentable-cats}. The author
used~\cite{adamek-rosicky} ,~\cite{borceux} Chapter V
of~\cite{stenstrom} to prepare this appendix.

We end this section by giving proofs to a few lemmas will be used
again and again throughout the paper. The first two are very basic
and concern generators in Grothendieck categories. We will usually
use Lemmas~\ref{lemma-lifting implies onto} and~\ref{lemma-class
contains a generator} without explicit mention. The last lemma is
a lifting property that the author learned from Mark Hovey through
personal communication.

\begin{lemma}\label{lemma-lifting implies onto}
Let $G$ be an object in an abelian category $\cat{G}$. $G$ is a
generator for $\cat{G}$ if and only if  given any morphism $d
\mathcolon C \ar D$, $d$ is an epimorphism whenever $d_*$ is an
epimorphism. Here $d_* \mathcolon \cat{G}(G,C) \ar \cat{G}(G,D)$
is defined by $d_*(t) = dt$.
\end{lemma}

\begin{proof}
Suppose $G$ is a generator and $d \mathcolon C \ar D$ is a
morphism for which $d_*$ is an epimorphism. To show $d$ is an
epimorphism, we will show that if $hd = 0$, then $h = 0$. By way
of contradiction suppose $h \neq 0$. By definition of a generator,
there exists a morphism $s \mathcolon G \ar D$ such that $hs \neq
0$. Now since $d_*$ is an epimorphism, there exists $t$ for which
$dt = s$. So now $0 = 0t = hdt = hs \neq 0$ which is a
contradiction. So $d$ must be an epimorphism.

Next, suppose $d$ is an epimorphism whenever $d_*$ is an
epimorphism and let $h \mathcolon X \ar Y$ be a nonzero map. To
show $G$ is a generator we need to find a morphism $s \mathcolon G
\ar X$ such that $hs \neq 0$. Let $k \mathcolon Z \ar X$ be the
kernel of $h$ and by way of contradiction suppose $hs = 0$ for all
$s$. By the universal property of kernel, each $s$ factors through
$k$. I.e., $k_*$ is an epimorphism. By hypothesis, $k$ is an
epimorphism. But since $k$ is a kernel it is also a monomorphism
and therefore $k$ is an isomorphism. This implies $h = 0$
\end{proof}

If $\cat{G}$ is an abelian category and $\class{A}$ is a class of
objects, we say that $\cat{G}$ has enough $\class{A}$-objects if
for any object $X \in \cat{G}$ we can find an epimorphism $A \ar
X$ where $A \in \class{A}$.

\begin{lemma}\label{lemma-class contains a generator}
Let $\cat{G}$ be a Grothendieck category and $\class{A}$ be a
class of objects which is closed under coproducts. Then
$\class{A}$ contains a generator if and only if $\cat{G}$ has
enough $\class{A}$-objects.
\end{lemma}

\begin{proof}
First assume $\class{A}$ contains a generator $G$ and that $C$ is
an arbitrary object in $\cat{C}$. It is well-known (for example,
see~\cite{borceux}) that the canonical map $$\bigoplus_{f \in
\Hom(G,C)} \hspace{-.1in} G \, \, \ar C$$ is an epimorphism. Since
we assume $\class{A}$ is closed under coproducts we are done.

Next assume $\cat{G}$ has enough $\class{A}$-objects and let $G$
be a generator for $\cat{G}$. Then we can find an epimorphism $A
\ar G$ for which $A \in \class{A}$. It is easy to check that $A$
too must be a generator for $\cat{G}$.
\end{proof}

\begin{lemma}\label{lemma-ext vanishing implies a lift}
Let $\cat{G}$ be any abelian category. Suppose we have short exact
sequences $A \xrightarrow{i} B \xrightarrow{p} C$ and $K
\xrightarrow{j} L \xrightarrow{q} M$ and a commutative diagram as
shown below:
$$\begin{CD}
      A @>f>>  L    \\
         @V i VV     @VVqV    \\
      B @>>g>  M \\
  \end{CD}$$
If $\, \Ext^1_{\cat{A}}(C,K) = 0$ then there exists a lift $h
\mathcolon B \ar L$ so that $hi = f$ and $qh = g$.
\end{lemma}

\begin{proof}
Consider the diagram below:
$$\begin{CD}
0       @>>>     A      @=      A \\
@VVV             @VV  \binom{i}{f} V           @VViV     \\
L       @> \binom{0}{1} >> B \oplus L @> (1 \ 0) >>    B \\
@VqVV            @VV (-g \ h) V           @VVV\\
M       @=       M      @>>>     0\\
\end{CD}$$
Each column forms a chain complex and so the diagram is a short
exact sequence of chain complexes. The associated long exact
sequence in homology leads to a short exact sequence  $K
\xrightarrow{k} T \xrightarrow{r} C$.

If we let $Z$ denote the kernel of the map $(-g \ q) \mathcolon B
\oplus L \ar M$, then $Z$ is the pullback of the maps $q$ and $g$
as shown in the square below:
$$\begin{CD}
      Z @>\tilde{g}>>  L    \\
         @V \tilde{q} VV     @VV q V    \\
      B @> g >>  M \\
  \end{CD}$$
The maps $\tilde{g}$ and $\tilde{q}$ are the projections $B \oplus
L \ar L$ and $B \oplus L \ar B$ restricted to $Z$.

Some more notation is necessary. We set $\tilde{\iota} :=
\binom{i}{f} \mathcolon A \ar Z$, and so $T =
\cok{\tilde{\iota}}$. We write $\tilde{p} \mathcolon Z \ar T$ for
the quotient map and we set $\tilde{k} := (0 \ j) \mathcolon K \ar
Z$. Now one can check that the diagram below commutes, the rows
and columns are exact, and the bottom right square is a pullback:
$$\begin{CD}
@.            K                   @=      K \\
@.           @V \tilde{k} VV              @VV k V     \\
A       @> \tilde{\iota} >> Z @> \tilde{p} >>  T \\
@|            @V \tilde{q} VV            @VV r V\\
A      @> i >>       B     @> p >>     C\\
\end{CD}$$

Now we can finally construct the lift. Since
$\Ext^1_{\cat{G}}(C,K) = 0$, the sequence $K \xrightarrow{k} T
\xrightarrow{r} C$ splits and so we have a map $n \mathcolon C \ar
T$ such that $rn = 1_C$. By the pullback property, there is a
unique map $\tilde{n} \mathcolon B \ar Z$ with $\tilde{q}\tilde{n}
= 1_B$ and $\tilde{p} \tilde{n} = np \mathcolon B \ar T$. We claim
that the map $h := \tilde{g} \tilde{n} \mathcolon B \ar L$ is a
lift. Indeed, we first have $qh = q \tilde{g} \tilde{n} = g
\tilde{q} \tilde{n} = g$. Next, one checks that $\tilde{p}
(\tilde{\iota} - \tilde{n} i) = 0$ and $\tilde{q}(\tilde{\iota} -
\tilde{n} i) = 0$ and so the pullback property tell us that
$\tilde{n} i = \tilde{\iota} \mathcolon A \ar Z$. This gives us
$hi = \tilde{g} \tilde{\iota} = f$ as required.
\end{proof}

\section{Small cotorsion pairs}\label{sec-small cotorsion pairs}

In the paper~\cite{hovey} we learned a relationship between
cotorsion pairs and model structures on an abelian category. In
that paper Hovey defined small cotorsion pairs, which are
cotorsion pairs $(\class{A},\class{B})$ along with a set of
``generating monomorphisms'' $I$. The basic analogy is that
\emph{small} is to cotorsion pairs as \emph{cofibrantly generated}
is to model categories.

\subsection{Properties of small cotorsion pairs}
The results in this subsection are all known and in fact much of
it can be found in Section~6 of~\cite{hovey}. The motivation for
this subsection is proving Lemmas~\ref{lemma-small cotorsion pair
properties} and~\ref{lemma-classification of left side of
cotorsion pair}. Although these were certainly known to Hovey
(personal communication) they do not appear in~\cite{hovey} and
the author could not find them in the literature.
Lemma~\ref{lemma-small cotorsion pair properties} will be used in
the next section while Lemma~\ref{lemma-classification of left
side of cotorsion pair} will be used to prove that the flat model
structure on chain complexes of quasi-coherent sheaves is
monoidal.

\begin{definition} Let $(\class{A},\class{B})$ be a
cotorsion pair in a Grothendieck category $\cat{G}$. Suppose
$\class{A}$ contains a generator $G$ for the category. The
following conditions are equivalent and we say that the cotorsion
pair is \emph{small} if it satisfies one of these conditions:

1) There is a set $\class{S}$ which cogenerates
$(\class{A},\class{B})$ and for each $S \in \class{S}$ there is a
monomorphism $i_S$, with $\cok{i_S} = S$, satisfying the
following: For all $X \in \cat{G}$, if $\cat{G}(i_S,X)$ is
surjective for all $S \in \class{S}$, then $X \in \class{B}$.

2) There is a set $I$ of monomorphisms for which $\cok{I} = \{
\cok{i} : i \in I \}$ cogenerates $(\class{A},\class{B})$ and
which satisfies the following: For all $X \in \cat{G}$, if
$\cat{G}(i,X)$ is surjective for all $i \in I$, then $X \in
\class{B}$.

3) There is a single monomorphism $i$ for which $\{ \cok{i} \}$
cogenerates $(\class{A},\class{B})$ and such that for all $X$,
$\cat{G}(i,X)$ surjective implies $X \in \class{B}$.

We call a collection $I$, as in (2) above, together with the
monomorphism $0 \ar G$, a set of \emph{generating monomorphisms}
for $(\class{A},\class{B})$.
\end{definition}

The first definition is the original due to Hovey and can be found
in~\cite{hovey}. The first definition clearly implies the second
definition. The second definition implies the third by looking at
the direct sum $\bigoplus_{i \in I} i$. Finally the third
condition clearly implies the first.

The following lemma can be found implicitly by studying Section~6
of~\cite{hovey}. See~\cite{hovey-model-categories} for the
definition of $I$-inj, $I$-cof, and $I$-cell.

\begin{lemma}\label{lemma-I-inj classification for small cot pairs}
Let $(\class{A},\class{B})$ be a cotorsion pair in a Grothendieck
category $\cat{G}$. Also suppose that $\class{A}$ contains a
generator $G$ for the category and that the cotorsion pair is
small with generating monomorphisms $I$. Then $I$-inj is the class
of all epimorphisms $p$ with $\ker{p} \in \class{B}$.
\end{lemma}

\begin{proof}
Say $p \mathcolon X \ar Y$ is in $I$-inj. Since $0 \ar G$ is in
$I$, there is a lift in any diagram of the form
$$\begin{CD}
0 @>>> X \\
@VVV  @VV p V \\
G @>>> Y \\
\end{CD}$$ By Lemma~\ref{lemma-lifting implies onto}, $p$
must be an epimorphism.

Now $\ker{p} \ar 0$ can be viewed as the pullback of $p$ over the
map $0 \ar Y$. Since $I$-inj is closed under pullbacks we see that
$\ker{p} \ar 0$ is also in $I$-inj. But this is equivalent to
saying $\cat{G}(i,\ker{p})$ is surjective for all $i \in I$. So
$\ker{p} \in \class{B}$ by the definition of a set of generating
monomorphisms.

On the other hand, let $p$ be an epimorphism with $\ker{p} \in
\class{B}$. We look for a lift in a diagram of the form
$$\begin{CD}
M @>>> X \\
@V i VV  @VV p V \\
N @>>> Y \\
\end{CD}$$ where $i \in I$. Since $\Ext(\cok{i},\ker{p}) = 0$,
such a lift exists Lemma~\ref{lemma-ext vanishing implies a lift}.
This proves the Lemma.
\end{proof}

\begin{lemma}\label{lemma-I-cof classification for small cot pairs}
Let $(\class{A},\class{B})$ be a cotorsion pair in a Grothendieck
category $\cat{G}$. Also suppose that $\class{A}$ contains a
generator $G$ for the category and that the cotorsion pair is
small with generating monomorphisms $I$. Then $I$-cof is the class
of all monomorphisms $f$ with $\cok{f} \in \class{A}$.
\end{lemma}

\begin{proof}
If $f$ is a monomorphism with $\cok{f} \in \class{A}$, then $f$ is
in $I$-cof by combining Lemma~\ref{lemma-I-inj classification for
small cot pairs} and Lemma~\ref{lemma-ext vanishing implies a
lift}.

Conversely, let $f \mathcolon M \ar N$ be in $I$-cof. Embed $M
\hookrightarrow I$ in an injective. Then $(I \ar 0) \in I$-inj, so
there is a lift in the diagram
$$\begin{CD}
M @>>> I \\
@V f VV @VVV \\
N @>>> 0 \\
\end{CD}$$ and so $f$ must be injective.

Now we wish to show that $N/M \in \class{A}$. So let $B \in
\class{B}$ be arbitrary and embed $B$ in an injective to get a
short exact sequence $0 \ar B \xrightarrow{i} I \xrightarrow{j}
I/B \ar 0$. Since $\Ext^1(N/M,I) = 0$, we will be done if we can
show that any map $h \mathcolon N/M \ar I/B$ lifts to a map
$\bar{h} \mathcolon N/M \ar I$ so that $h = j \bar{h}$. Notice we
have a commutative diagram
$$\begin{CD}
M @> 0 >> I @= I \\
@V f VV @. @VV j V \\
N @>> c > N/M @>> h > I/B \\
\end{CD}$$ where $c \mathcolon N \ar N/M$ is the canonical map to the cokernel.
Since $j \in I$-inj and $f \in I$-cof, there exists a lift $\psi
\mathcolon N \ar I$ such that $\psi f = 0$ and $j \psi = h c$. By
the universal property of cokernel, there exists a map $\bar{h}
\mathcolon N/M \ar I$ such that $\psi = \bar{h} c$. But now $hc =
j \psi = j \bar{h} c$, and since $c$ is epi, it is right
cancellable. Thus $h = j \bar{h}$.
\end{proof}

The next lemma basically says that a set $I$ of monomorphisms (for
which $\cok{I} = \class{S}$ cogenerates $(\class{A},\class{B})$)
is a set of generating monomorphisms if and only if $I$-inj and
$I$-cof can be classified as in the last two lemmas.

\begin{lemma}\label{lemma-small cotorsion pair properties}
Let $(\class{A},\class{B})$ be a cotorsion pair in a Grothendieck
category $\cat{G}$, cogenerated by a set $\class{S}$. Also suppose
that $\class{A}$ contains a generator $G$ for the category.
Suppose for each $S \in \class{S}$ there is a monomorphism $i_S$,
with $\cok{i_S} = S$. Denote the set of all $i_S$ together with $0
\ar G$ by $I$. Then the following are equivalent:

1) $(\class{A},\class{B})$ is a small cotorsion pair with $I$ a
set of generating monomorphisms.

2) $I$-inj is the class of all epimorphisms with kernel in
$\class{B}$.

3) $I$-cof is the class of all monomorphisms with cokernel in
$\class{A}$.

\end{lemma}

\begin{proof}
(1) implies (2) is Lemma~\ref{lemma-I-inj classification for small
cot pairs} and (1) implies (3) is Lemma~\ref{lemma-I-cof
classification for small cot pairs}. We first show (2) implies
(1). So let $X \in \cat{G}$ and suppose $\cat{G}(i_S,X)$ is
surjective for all $S \in \class{S}$. This is equivalent to saying
$X \ar 0$ is in $I$-inj. By hypothesis $X \in \class{B}$ as
desired.

Next we show (3) implies (1). So let $X \in \cat{G}$ and suppose
$\cat{G}(i_S,X)$ is surjective for all $S \in \class{S}$. Again,
this is equivalent to saying $X \ar 0$ is in $I$-inj. We need to
conclude $X \in \class{B}$, so we let $A \in \class{A}$ be
arbitrary and argue that $$0 \ar X \xrightarrow{f} Z
\xrightarrow{g} A \ar 0$$ must split. But by hypothesis, $f \in
I$-cof = $I$-inj-proj, so there is a lift in the diagram
$$\begin{CD}
X @= X \\
@V f VV  @VVV \\
Z @>>> 0 \\
\end{CD}$$
This lift is a splitting.

\end{proof}

\begin{lemma}\label{lemma-classification of left side of cotorsion pair}
Let $(\class{A},\class{B})$ be a cotorsion pair in a Grothendieck
category $\cat{G}$. Also suppose that $\class{A}$ contains a
generator $G$ for the category and that the cotorsion pair is
small with cogenerating set $\class{S}$ and generating
monomorphisms $I$. Then every object $A \in \class{A}$ is a
retract of a transfinite extension of objects in $\class{S}$. In
particular, $\class{A}$ is the smallest class containing
$\class{S}$ and closed under transfinite extensions and retracts
(summands).
\end{lemma}

\begin{proof}
Let $A \in \class{A}$. Then $0 \ar A$ is in $I$-cof by
Lemma~\ref{lemma-I-cof classification for small cot pairs}. By
Corollary~2.1.15 of~\cite{hovey-model-categories} we see $0 \ar A$
is a retract of a map $0 \ar Y$ in $I$-cell, by a map which fixes
0. That is, $A$ is a retract of $Y$. But if $0 \ar Y$ is in
$I$-cell, then $Y$ is a transfinite extension of the cokernels of
maps in $I$. (This follows right from the definition of $I$-cell
and the fact that cokernels are unchanged when pushing out over a
monomorphism.) So $A$ is a retract of a transfinite extension of
objects in $\class{S}$.

Say $\class{W}$ is the smallest class containing $\class{S}$ and
is closed under transfinite extensions and retracts. Now
$\class{A}$ contains $\class{S}$ and we know from Lemma~6.2
of~\cite{hovey} that the left side of a cotorsion pair is always
closed under transfinite extensions and retracts. So $\class{A}
\supseteq \class{W}$. Conversely, from the last paragraph it is
clear that $\class{A} \subseteq \class{W}$.
\end{proof}

\subsection{An associated cotorsion pair of chain complexes is
small}\label{subsec-an associated cotorsion pair of chain
complexes is small}

We continue to let $(\class{A},\class{B})$ represent a small
cotorsion pair in a Grothendieck category $\cat{G}$ which has a
generator $G \in \class{A}$. By Corollary~3.8 of~\cite{gillespie},
we have induced cotorsion pairs $(\ldgclass{A},\rclass{B})$ and
$(\lclass{A},\rdgclass{B})$ of chain complexes. We now show
$(\ldgclass{A},\rclass{B})$ is small whenever
$(\class{A},\class{B})$ is small. The author doesn't see a
corresponding theorem for the cotorsion pair
$(\lclass{A},\rdgclass{B})$ without making further assumptions on
the class $\class{A}$. This is the subject of
Section~\ref{sec-locally cogenerated and kaplansky}.

\begin{lemma}\label{lemma-lifting implies exactness}
Let $X$ be a chain complex in an abelian category $\class{G}$ with
generator $G$. If any chain map $f \mathcolon S^n(G) \ar X$
extends to $D^{n+1} (G)$, then $X$ is exact.
\end{lemma}

\begin{proof}
Let $n$ be an arbitrary integer. By Lemma~\ref{lemma-lifting
implies onto}, showing exactness in degree $n$ requires showing
that any morphism $f \mathcolon G \ar Z_nX$ lifts over $d
\mathcolon X_{n+1} \ar Z_nX$. But it is easy to see that this is
the same as showing that the induced chain map $\hat{f} \mathcolon
S^n(G) \ar X$ extends to a morphism $D^{n+1}(G) \ar X$.
\end{proof}

\begin{proposition}\label{prop-induced cotorsion pair of complexes is small}
Let $(\class{A},\class{B})$ be a cotorsion pair in a Grothendieck
category $\cat{G}$ which has a generator $G \in \class{A}$. If
$(\class{A},\class{B})$ is cogenerated by a set $\{A_i\}_{i \in
I}$, then the induced cotorsion pair $(\ldgclass{A},\rclass{B})$
is cogenerated by the set
$$S = \{\, S^n(G) \, | \, n \in \Z \,\} \cup \{\, S^n(A_i) \, | \, n \in \Z \, , \, i \in I \,
\}.$$

Furthermore, suppose $(\class{A},\class{B})$ is small with
generating monomorphisms the map $0 \ar G$ together with
monomorphisms $k_i$ as below (one for each $i \in I$): $$0 \ar Y_i
\xrightarrow{k_i} Z_i \ar A_i \ar 0.$$ Then
$(\ldgclass{A},\rclass{B})$ is small with generating monomorphisms
the set $$I = \{\, 0 \ar D^n(G) \, \} \cup \{\, S^{n-1}(G) \ar
D^n(G) \,\} \cup \{\, S^n(Y_i) \xrightarrow{S^n(k_i)} S^n(Z_i) \,
\}.$$
\end{proposition}

\begin{proof}
Clearly $S \subseteq \ldgclass{A}$, so we have $\rightperp{S}
\supseteq \rightperp{(\ldgclass{A})} = \rclass{B}$. Conversely if
$X \in \rightperp{S}$, then $0 = \Ext^1_{\text{Ch}(\cat{G})}
(S^n(A_i),X)$ for all $i \in I$. But $\Ext^1_{\text{Ch}(\cat{G})}
(S^n(A_i),X) \cong \Ext^1_{\cat{G}} (A_i,Z_nX)$ (Lemma~3.1
in~\cite{gillespie}). So $\Ext^1_{\cat{G}} (A_i,Z_nX) = 0$ which
implies $Z_n X \in \class{B}$ since the set $\{A_i\}$ cogenerates
the cotorsion theory.

Next we want to show that $X$ is exact. Consider the short exact
sequence
$$0 \ar S^{n-1}(G) \ar D^{n}(G) \ar S^{n}(G) \ar 0.$$ It induces an exact sequence
of abelian groups $$\Hom_{\text{Ch}(\cat{G})} (D^{n}(G),X) \ar
\Hom_{\text{Ch}(\cat{G})} (S^{n-1}(G),X) \ar
\Ext_{\text{Ch}(\cat{G})} (S^{n}(G),X).$$ But again Lemma~3.1
of~\cite{gillespie} gives us $\Ext_{\text{Ch}(\cat{G})}
(S^{n}(G),X) \cong \Ext_{\cat{G}} (G,Z_{n} X)$ and this last group
equals 0 by the last paragraph. Therefore Lemma~\ref{lemma-lifting
implies exactness} tells us $X$ is exact. Since $X$ is exact and
has cycles in $\class{B}$, we see that $X \in \rclass{B}$. So
$\rightperp{S} = \rclass{B}$. This shows that $\class{S}$
cogenerates the cotorsion pair $(\ldgclass{A},\rclass{B})$.

Next we prove the statement about smallness. First note that since
$G$ generates $\cat{G}$, the complexes $D^n(G)$ generate
Ch($\cat{G}$). Also $D^n(G) \in \ldgclass{A}$, and so
$\ldgclass{A}$ contains the generators $\{ D^n(G) \}$. Now let $X$
be any chain complex. We wish to show that ``extending through
monomorphisms in $I$'' implies $X \in \rclass{B}$. But again, if
morphisms $S^{n-1}(G) \ar X$ can extend over $D^n(G)$, then $X$
must be exact by Lemma~\ref{lemma-lifting implies exactness}. Next
let $Z_nX$ be a cycle. Any map $Y_i \ar Z_n X$ determines a
morphism $S^n(Y_i) \ar X$, which we assume extends over $S^n(k_i)$
to a map $S^n(Z_i) \ar X$. Thus any map $Y_i \ar Z_n X$ extends
over $k_i$ to a map $Z_i \ar Z_n X$. By hypothesis this implies
$Z_nX \in \class{B}$.
\end{proof}

As mentioned at the beginning of this subsection, we know from
Corollary~3.8 of~\cite{gillespie} that if $\class{A}$ contains a
generator for $\cat{G}$ and $(\class{A},\class{B})$ is a cotorsion
pair, then we have two induced cotorsion pairs on Ch($\class{A}$).
We denote them by $(\ldgclass{A},\rclass{B})$ and
$(\lclass{A},\rdgclass{B})$ as in~\cite{gillespie}. We will see in
the proof of Theorem~\ref{them-model structure} that a crucial
step in building a model category structure on Ch($\class{A}$) is
showing that these induced cotorsion pairs are \emph{compatible}.
This means $\ldgclass{A} \cap \class{E} = \lclass{A}$ and
$\rdgclass{B} \cap \class{E} = \rclass{B}$. Condition (4) of
Corollary~\ref{cor-equivalence of compatible and coresolving}
below will be used at that time to guarantee that the induced
cotorsion pairs are compatible.

\begin{corollary}\label{cor-equivalence of compatible and
coresolving} Let $(\class{A},\class{B})$ be a small cotorsion pair
in a Grothendieck category $\cat{G}$. Assume $\class{A}$ contains
a generator for $\cat{G}$. Then the following are equivalent:

(1) The induced cotorsion pairs $(\ldgclass{A},\rclass{B})$ and
$(\lclass{A},\rdgclass{B})$ are compatible.

(2) $\class{A}$ is resolving and $\class{B}$ is coresolving. That
is, $\class{A}$ is closed under taking kernels of epimorphisms and
$\class{B}$ is closed under taking cokernels of monomorphisms.

(3) $\Ext^n_{\class{A}} (A,B) = 0$ for any $n>0$ and any $A \in
\class{A}$ and $B \in \class{B}$.

(4) $\lclass{A} = \ldgclass{A} \cap \class{E}$
\end{corollary}

\begin{proof}
First we show (1) implies (2). Say we are given a short exact
sequence $0 \ar X \ar A_0 \ar A \ar 0$ where $A_0$ and $A$ are in
the class $\class{A}$. Since $\class{A}$ contains a generator we
can complete the short exact sequence to obtain an
``$\class{A}$-resolution'' $\cdots \ar A_2 \ar A_1 \ar A_0 \ar A
\ar 0$. Since this resolution is a bounded below complex with
objects in $\class{A}$, it is in $\ldgclass{A}$ by Lemma~3.4
of~\cite{gillespie}. But it is also in $\class{E}$. Since
$\ldgclass{A} \cap \class{E} = \lclass{A}$, we see that the
resolution is in $\lclass{A}$. Therefore $X \in \class{A}$ by the
definition of the class $\lclass{A}$. This shows $\class{A}$ is
resolving. The dual shows that $\class{B}$ is coresolving.

Next we show (2) implies (3). We will prove by induction that
$\Ext^n(A,B) = 0$ for all $n
> 0$ , $A \in \class{A}$ and $B \in \class{B}$.
Obviously $\Ext^1(A,B) = 0$ for any $A \in \class{A}$ and
$\class{B}$. Now suppose $k > 0$ and that $\Ext^k(A,B) = 0$ for
any $A \in \class{A}$ and $B \in \class{B}$. We now let $A \in
\class{A}$ and $B \in \class{B}$ be arbitrary but fixed and wish
to argue $\Ext^{k+1}(A,B) = 0$. We start by embedding $B$ inside
an injective, $0 \ar B \ar I \ar B' \ar 0$. Note that $I \in
\class{B}$ and so by hypothesis $B' \in \class{B}$. We now apply
the functor $\Hom(A,-)$ to get the long exact sequence $\cdots \ar
\Ext^k(A,B') \ar \Ext^{k+1}(A,B) \ar \Ext^{k+1}(A,I) \ar \cdots$.
By the induction hypothesis $\Ext^k(A,B') = 0$ and since $I$ is
injective $\Ext^{k+1}(A,I) = 0$. It follows that $\Ext^{k+1}(A,B)
= 0$. Therefore (2) implies (3).

It is easy to see that (3) implies (2), so we have that (2) and
(3) are equivalent.

Now we show (2) implies (4). First use Lemma~3.10
of~\cite{gillespie} to see that $\lclass{A} \subset \ldgclass{A}
\cap \class{E}$. Next, using only the coresolving hypothesis in
(2), we can perform the (dual of the) argument in the proof of
Theorem~3.12 of~\cite{gillespie} to conclude $\lclass{A} \supset
\ldgclass{A} \cap \class{E}$. So $\lclass{A} = \ldgclass{A} \cap
\class{E}$.

It is left to show that (4) implies (1) and this amounts to
showing $\rclass{B} = \rdgclass{B} \cap \class{E}$. But by
Proposition~\ref{prop-induced cotorsion pair of complexes is
small} we see that $(\ldgclass{A} , \rclass{B})$ is a small
cotorsion pair. By Corollary~6.6 of~\cite{hovey} we get that
$(\ldgclass{A},\rclass{B})$ has enough injectives. Finally
Lemma~3.14 of~\cite{gillespie} tells us $\rclass{B} = \rdgclass{B}
\cap \class{E}$.
\end{proof}


\section{Locally cogenerated classes and Kaplansky
classes}\label{sec-locally cogenerated and kaplansky}

Throughout this section we again assume that $\cat{G}$ is a
Grothendieck category. In each of the categories $R$-Mod,
$\cat{O}_X$-Mod and $\qcox$, the class of flat objects satisfies
an important property. In a loose sense, every flat object is
``built'' from ``smaller'' flat objects. ''Built'' in this case
means ``is a transfinite extension of'' and by ``smaller'' we
essentially mean ``of smaller cardinality''. We axiomatize this to
get the concept of a locally cogenerated class. We will prove that
a locally cogenerated class $\class{F}$ which is closed under
transfinite extensions and retracts gives rise to a small
cotorsion pair $(\class{F},\class{C})$. However, in light of
Section~\ref{subsec-an associated cotorsion pair of chain
complexes is small}, to get the flat model structure (on either
$R$-Mod, $\cat{O}_X$-Mod or $\qcox$) we really need the induced
cotorsion pair $(\lclass{F}, \rdgclass{C})$ to be small. We will
do this by showing $\lclass{F}$ is a locally cogenerated class.
Doing so will lead us to the notion of a Kaplansky class, which is
a slight strengthening of a locally cogenerated class. The author
got the term ``Kaplansky class'' from Edgar Enochs whose work
(with several coauthors) on proving the existence of flat covers
in categories such as $R$-Mod, $\cat{O}_X$-Mod and $\qcox$ is
closely related to the idea of a Kaplansky class.
See~\cite{enochs-flat-cover-theorem}~,
\cite{enochs-kaplansky-classes}~, \cite{enochs-quasi-coherent}~,
\cite{enochs-quasi-coherent-projectiveline}~,
and~\cite{enochs-oyonarte}. In each category mentioned above,
Enochs and coauthors have essentially shown that the class of flat
objects form a Kaplansky class. So as a result of
Theorem~\ref{them-model structure} we then have an induced (flat)
model structure on the associated chain complex category. We
explain in section~\ref{sec-quasi-coherent sheaves} using results
of~\cite{enochs-quasi-coherent} why the class of flat objects in
$\qcox$ is a Kaplansky class.

The reader may want to skim through
Appendix~\ref{appendix-locally-presentable-cats} at this point
since many definitions and theorems we use in this section can be
found there.

\subsection{Locally cogenerated classes}\label{subsec-locally cogenerated classes}

\begin{definition}
Let $\kappa$ be a regular cardinal. Given a class $\class{F}$ of
objects in $\cat{G}$, we say $\class{F}$ is \emph{locally
$\kappa$-cogenerated} if for every $0 \neq F \in \class{F}$, there
exists $0 \neq S \subseteq F$ with $S \in \mathbf{Gen}_{\kappa} \
\class{F}$ and $F/S \in \class{F}$. We say $\class{F}$ is
\emph{locally cogenerated} if it is locally $\kappa$-cogenerated
for some regular cardinal $\kappa$.
\end{definition}

\begin{fact}\label{fact-e}
Let $\kappa'$ and $\kappa$ be regular cardinals with $\kappa' \geq
\kappa$. It follows from Fact~\ref{fact-a} that if $\class{F}$ is
locally $\kappa$-cogenerated, then $\class{F}$ is locally
$\kappa'$-cogenerated.
\end{fact}

\begin{lemma}\label{lemma-write B as a transfinite ext} Suppose
the class $\class{F}$ is locally $\kappa$-cogenerated. Then given
a monomorphism $f \mathcolon A \hookrightarrow B$ with $B/A \in
\class{F}$, we may write $B$ as a transfinite extension of $A$ by
objects in $\mathbf{Gen}_{\kappa} \ \class{F}$. In particular,
every object of $\class{F}$ is a transfinite extension of
$\mathbf{Gen}_{\kappa} \ \class{F}$
\end{lemma}

\begin{proof}
 Set $X_0 = A$. Now use the
definition of locally $\kappa$-cogenerated to find $0 \neq X_1/X_0
\subseteq B/X_0$ with $X_1/X_0  \in \mathbf{Gen}_{\kappa} \
\class{F}$ and $B/X_1 \in \class{F}$. Continue, by transfinite
induction, setting $X_{\gamma} = \bigcup_{\alpha < \gamma}
X_{\alpha}$ for limit ordinals $\gamma$. Eventually, this process
must terminate (since $B$ has only a set of subobjects by
Fact~\ref{fact-wellpowered}), and we get $B = \bigcup_{\alpha <
\lambda} X_{\alpha}$ for some ordinal $\lambda$. This says $B$ is
a transfinite extension of $A$ by objects in
$\mathbf{Gen}_{\kappa} \ \class{F}$ since $X_0 = A$ and
$X_{\alpha+1}/X_{\alpha}  \in \mathbf{Gen}_{\kappa} \ \class{F}$.
\end{proof}

The author learned the following Lemma from Mark Hovey. It is a
generalization of Lemma~V.3.3 from~\cite{stenstrom}.

\begin{lemma}\label{surjective-lemma}
Suppose $\cat{G}$ is locally $\kappa$-generated. Given an
epimorphism $g \mathcolon X \ar Y$ where $Y$ is
$\kappa$-generated, there exists a $\kappa$-generated subobject
$X' \subseteq X$ for which $g_{|X'} \mathcolon X' \ar Y$ is an
epimorphism.
\end{lemma}

\begin{proof}
Since $\cat{G}$ is locally $\kappa$-generated, We may write $X =
\sum_{i \in I}X_i$ as a $\kappa$-filtered union of
$\kappa$-generated subobjects of $X$. Since $g$ is an epimorphism,
$Y = \sum_{i \in I} g(X_i)$, and this too is a $\kappa$-filtered
union. Now we must have $Y = g(X_i)$ for some $i \in I$ by
Fact~\ref{kappa-generated-characterization}. So $g_{|X_i}
\mathcolon X_i \ar Y$ is an epimorphism.
\end{proof}

\begin{definition}\label{def-the set I}
Suppose $\class{F}$ is locally cogenerated class in $\cat{G}$. By
Facts~\ref{fact-c.5} and~\ref{fact-e} we can choose a regular
cardinal $\kappa$ with each of the following properties:
 (1) $\cat{G}$ is locally $\kappa$-presentable (hence locally
$\kappa$-generated too) and (2) $\class{F}$ is locally
$\kappa$-cogenerated. Having chosen such a $\kappa$, we define $I$
to be the set of all (representatives of isomorphism classes of)
monomorphisms $A \hookrightarrow B$ for which $B \in
\mathbf{Gen}_{\kappa} \ \cat{G}$  and $B/A \in \class{F}$.
Furthermore, if $\class{F}$ contains a generator $G$, then we also
assume we have chosen $\kappa$ large enough so that (3) $G$ is
$\kappa$-presentable.
\end{definition}

Note that by Fact~\ref{fact-b.5}, if $(A \hookrightarrow B) \in I$
 then we have $B/A \in \mathbf{Gen}_{\kappa} \ \cat{F}$. Also, $(0
 \hookrightarrow B) \in I$ iff $B \in \mathbf{Gen}_{\kappa} \
\cat{F}$.

\begin{lemma}\label{characterization of I-cell}
Let $\class{F}$ be a locally cogenerated class in $\cat{G}$. If $A
\hookrightarrow B$ is a monomorphism with $B/A \in \class{F}$ then
$f \mathcolon A \ar B \in I$-cell. The converse holds if we assume
that $\class{F}$ is closed under transfinite extensions.
\end{lemma}

\begin{proof}
 We first need to show that any $f \mathcolon A
\hookrightarrow B$ with $B/A \in \class{F}$ is in $I$-cell. Using
Lemma~\ref{lemma-write B as a transfinite ext}, we may write $B =
\bigcup_{\alpha < \lambda} X_{\alpha}$ in such a way that $X_0 =
A$ and $X_{\alpha+1}/X_{\alpha} \in \mathbf{Gen}_{\kappa} \
\class{F}$. So $f$ is the composition of the $\lambda$-sequence
$$X_0 = A \hookrightarrow X_1
\hookrightarrow X_2 \cdots \xrightarrow{} X_{\alpha}
\hookrightarrow X_{\alpha+1} \xrightarrow{} \cdots$$ where each
$X_{\alpha} \hookrightarrow X_{\alpha+1}$ satisfies
$X_{\alpha+1}/X_{\alpha} \in \mathbf{Gen}_{\kappa} \ \class{F}$.
But since $X_{\alpha+1}/X_{\alpha}$ is $\kappa$-generated we can,
by Lemma~\ref{surjective-lemma}, find $X'_{\alpha+1} \subseteq
X_{\alpha+1}$ such that $X'_{\alpha+1}$ is $\kappa$-generated and
$X'_{\alpha+1} \ar X_{\alpha+1}/X_{\alpha}$ is surjective. Then
note $X_{\alpha+1} = X'_{\alpha+1} + X_{\alpha}$. Now we have the
following commutative diagram, where the rows are short exact
sequences and the left square is a pushout (and pullback):

$$\begin{CD}
0 @>>> X_{\alpha} @>>> X_{\alpha} + X'_{\alpha+1} @>>> X_{\alpha+1}/X_{\alpha} @>>> 0 \\
@. @AAA @AAA @| @.\\
0 @>>> X_{\alpha} \cap X'_{\alpha+1} @>>> X'_{\alpha+1} @>>> X_{\alpha+1}/X_{\alpha} @>>> 0 \\
\end{CD}$$
Now $X_{\alpha} \cap X'_{\alpha+1} \hookrightarrow X'_{\alpha+1}$
is in $I$ by definition. This shows $f$ is in $I$-cell, since $f$
is a transfinite composition of pushouts of maps in $I$.

To prove the converse, we assume that $\class{F}$ is closed under
transfinite extensions. (However, we do not need to assume that
$\class{F}$ is locally cogenerated for this direction.) Indeed
suppose $f \mathcolon A \ar B$ is in $I$-cell, so that it is the
composition of some $\lambda$-sequence $A = X_0 \ar \cdots
X_{\alpha} \ar X_{\alpha +1} \cdots$ where each $X_{\alpha} \ar
X_{\alpha +1}$ is a pushout of some $A_{\alpha} \hookrightarrow
A_{\alpha+1}$ in $I$. Since a pushout of a mono is a mono, each
$X_{\alpha} \hookrightarrow X_{\alpha+1}$ must be a monomorphism.
It follows that $f$, being the map from $A$ to the colimit of the
$\lambda$-sequence, $B$, is also a monomorphism.

Now take the quotient of the entire diagram, by $A = X_0$. It
gives a new cone whose base is the $\lambda$-sequence $0
\hookrightarrow X_1/X_0 \hookrightarrow \cdots X_{\alpha}/X_0
\hookrightarrow X_{\alpha+1}/X_0 \cdots$. All maps in the new
diagram are also monomorphisms and the universal property of a
quotient object will show that the new cone is actually a colimit
cone, with colimit $B/A$. Furthermore the new diagram shows that
$B/A$ is a transfinite extension of things in $\class{F}$, since
each $(X_{\alpha+1}/X_0) / (X_{\alpha}/X_0) \cong
X_{\alpha+1}/X_{\alpha}$ is in the class $\class{F}$ as the
following pushout diagram shows:

$$\begin{CD}
0 @>>> X_{\alpha} @>>> X_{\alpha+1} @>>> A_{\alpha+1}/A_{\alpha} @>>> 0 \\
@. @AAA @AAA @| @.\\
0 @>>> A_{\alpha} @>>> A_{\alpha+1} @>>> A_{\alpha+1}/A_{\alpha} @>>> 0 \\
\end{CD}$$
\end{proof}

\begin{remark}\label{remark on I-cell}
Let $\class{F}$ be a locally cogenerated class in $\cat{G}$. If $F
\in \class{F}$ then $0 \ar F \in I$-cell. The converse holds if
$\class{F}$ is closed under transfinite compositions.
\end{remark}

\begin{lemma}\label{I-inj classification}
Let $\class{F}$ be locally cogenerated and contain a generator
$G$. Then $I$-inj equals the class of all maps $p \mathcolon X \ar
Y$ such that $p$ is surjective and $\ker{p} \in
\rightperp{(\mathbf{Gen}_{\kappa} \ \class{F})}$
\end{lemma}

\begin{proof}
Suppose $p \mathcolon X \ar Y$ is in $I$-inj. First we show that
$p$ must be surjective. Using Lemma~\ref{lemma-lifting implies
onto} this is easy: $0 \ar G \in I$ and so there is a lift in the
diagram below.
$$\begin{CD}
0 @>>> X \\
@VVV @VVpV \\
G @>>> Y \\
\end{CD}$$ Next we show $\ker{p} \in \rightperp{(\mathbf{Gen}_{\kappa} \
\class{F})}$. Let $F \in \mathbf{Gen}_{\kappa} \ \class{F}$ be
arbitrary. We want to show that any extension $0 \ar \ker{p}
\xrightarrow{f} Z \xrightarrow{g} F \ar 0$ must split. But $f \in
I$-cell $\subseteq I$-cof, so there exists a lift in the
commutative diagram
$$\begin{CD}
\ker{p} @>>> X \\
@V f VV @VV p V \\
Z @>> 0 > Y. \\
\end{CD}$$ Call the lift $h$. Since $ph = 0$, $h$ lands in
$\ker{p}$ and so provides a retraction for $f$. This shows
$\ker{p} \in \rightperp{(\mathbf{Gen}_{\kappa} \ \class{F})}$.

Conversely, say $p \mathcolon X \ar Y$ is a surjection and
$\ker{p} \in \rightperp{(\mathbf{Gen}_{\kappa} \ \class{F})}$.
Given any map $f \mathcolon A \ar B $ in $I$, and a commutative
diagram
$$\begin{CD}
A @>>> X \\
@V f VV @VV p V \\
B @>>> Y \\
\end{CD}$$ we seek a lift. But since $\Ext^1(B/A, \ker{p}) = 0$,
Lemma~\ref{lemma-ext vanishing implies a lift} provides the lift!
\end{proof}

\begin{proposition}\label{prop-locally cogeneratd class induces small cotorsion pair}
Suppose $\class{F}$ is a locally cogenerated class containing a
generator $G$. Furthermore, suppose $\class{F}$ is closed under
transfinite extensions and retracts. Then $(\class{F},
\rightperp{\class{F}})$ is a small cotorsion pair with $I$ as the
set of generating monomorphisms.
\end{proposition}

\begin{proof}
First we show $(\class{F}, \rightperp{\class{F}})$ is a cotorsion
pair cogenerated by the set $\mathbf{Gen}_{\kappa} \ \class{F}$.
So we show $\class{F} =
\leftperp{(\rightperp{(\mathbf{Gen}_{\kappa} \ \class{F})})}$.
First let $F \in \class{F}$. By Lemma~\ref{lemma-write B as a
transfinite ext} $F$ is a transfinite extension of objects in
$\mathbf{Gen}_{\kappa} \ \class{F}$. It follows from the proof of
Lemma~6.2 in~\cite{hovey} that $\Ext(F,X) = 0$ for any $X \in
\rightperp{(\mathbf{Gen}_{\kappa} \ \class{F})}$. So $F \in
\leftperp{(\rightperp{(\mathbf{Gen}_{\kappa} \ \class{F})})}$.
Conversely, say $X \in
\leftperp{(\rightperp{(\mathbf{Gen}_{\kappa} \ \class{F})})}$.
Factor $0 \ar X$ using the small object argument as $0
\xrightarrow{i}  F \xrightarrow{p} X$ where $i \in I$-cell and $p
\in I$-inj. It follows from Remark~\ref{remark on I-cell} that $F
\in \class{F}$. It follows from Lemma~\ref{I-inj classification}
that $\ker{p} \in \rightperp{(\mathbf{Gen}_{\kappa} \
\class{F})}$. Since $\Ext(X,\ker{p}) = 0$ we see $p$ is a split
epimorphism. Therefore $X$ is a retract of $F$. By the retract
hypothesis $X \in \class{F}$.

Now it follows from Lemma~\ref{I-inj classification} and
Lemma~\ref{lemma-small cotorsion pair properties} that
$(\class{F},\rightperp{\class{F}})$ is a small cotorsion pair with
generating monomorphisms the set $I$.
\end{proof}


\subsection{Kaplansky classes}\label{subsec-kaplansky classes}

We now strengthen the idea of a locally cogenerated class to get
the idea of a Kaplansky class. Our goal then is to show that a
Kaplansky class in $\cat{G}$ gives rise to a model structure on
$\textnormal{Ch}(\cat{G})$. The term ``Kaplansky class'' first
appeared in~\cite{enochs-kaplansky-classes}. There it is defined
for modules over a ring $R$ as a class $\class{F}$ for which there
is a regular cardinal $\kappa$ with the following property: Given
$F \in \class{F}$ and $x \in F$, there exists an $S \in \class{F}$
with $x \in S \subseteq F$ and $|S| \leq \kappa$ and $F/S \in
\class{F}$. Below we formulate a version for abelian categories
and then prove Theorem~\ref{them-model structure}, the existence
of an induced model category structure.

\begin{definition}\label{def-kaplansky class}
Let $\class{F}$ be a class of objects in an abelian category and
let $\kappa$ be a regular cardinal. We say $\class{F}$ is a
\emph{$\kappa$-Kaplansky class} if the following property holds:
Given $X \subseteq F \neq 0$ where $F \in \class{F}$ and $X$ is
$\kappa$-generated, there exists $0 \neq S \in \class{F}$ such
that $S$ is $\kappa$-presentable and $X \subseteq S \subseteq F$
and $F/S \in \class{F}$. We say $\class{F}$ is a \emph{Kaplansky
class} if it is a $\kappa$-Kaplansky class for some regular
cardinal $\kappa$.
\end{definition}

The requirement that we can find a $\kappa$-presentable object
containing a $\kappa$-generated object in the definition of a
$\kappa$-Kaplansky class may seem strange. However, in practice
one usually needs to choose $\kappa$ to be very large in order to
show that $\class{F}$ is a $\kappa$-Kaplansky class. So why not
just take $\kappa$ to be so large that the $\kappa$-generated and
$\kappa$-presentable objects coincide? By
Appendix~\ref{appendix-kappa-generated and kappa-presentable
objects coincide} we can always find such a $\kappa$ as long as we
are in a Grothendieck category. When this is the case,
$\kappa$-generated becomes a good substitute for the notion of
cardinality. For example, subobjects and quotient objects of
$\kappa$-generated objects will again be $\kappa$-generated. (This
follows from facts in
Appendix~\ref{appendix-locally-presentable-cats}.) The method we
just outlined is exactly what we will use in
Section~\ref{sec-quasi-coherent sheaves} to show that the class of
flat quasi-coherent sheaves form a Kaplansky class. At any rate,
we will need the definition of a $\kappa$-Kaplansky class to be as
stated above in order to prove Proposition~\ref{prop-F-twiddle
locally cogenerated}.

It is obvious that if $\class{F}$ is a $\kappa$-Kaplansky class,
then $\class{F}$ is locally $\kappa$-cogenerated. We next prove
that if $\class{F}$ is a Kaplansky class, then the class of
$\class{F}$-complexes, denoted $\lclass{F}$~\cite{gillespie},
ought to be locally cogenerated as well.

\begin{lemma}\label{generated complexes}
Let $\cat{G}$ be an abelian category and let $X$ be a chain
complex in $\textnormal{Ch}(\cat{G})$. If $\kappa > \omega$ is a
regular cardinal, then $X$ is $\kappa$-generated if and only if
$X_n$ is $\kappa$-generated for each $n$. $X$ is
$\omega$-generated (finitely generated) if and only if $X$ is
bounded, above and below, and each $X_n$ is $\omega$-generated
(finitely generated).
\end{lemma}

\begin{proof}
Say $\kappa$ is \emph{any} regular cardinal, and let $n$ be an
arbitrary integer. Let $X_n = \sum_{i \in I} X_i$ where $X_i
\subseteq X_n$. By Fact~\ref{kappa-generated-characterization} we
wish to show $X_n = \sum_{i \in J} X_i$ where $J \subseteq I$ and
$|J| < \kappa$. Each $X_i$ gives rise to a subcomplex $C(X_i)$ of
$X$ by defining $$C(X_i) = \cdots \ar
d_{n+2}^{-1}(d_{n+1}^{-1}(X_i)) \ar d_{n+1}^{-1}(X_i) \ar X_i \ar
X_{n-1} \ar X_{n-2} \ar \cdots$$ and furthermore $X = \sum_{i \in
I} C(X_i)$. Since $X$ is $\kappa$-generated,
Fact~\ref{kappa-generated-characterization} tells us $X = \sum_{i
\in J} C(X_i)$ where $J \subseteq I$ and $|J| < \kappa$. Thus $X_n
= \sum_{i \in J} X_i$. This shows that if $X$ is
$\kappa$-generated, then $X_n$ is $\kappa$-generated in $\cat{G}$.
For the special case of when $\kappa = \omega$, we still need to
argue that $X$ is bounded. For this, write $X = \sum_{n = \Z} S_n$
where $$S_n = \cdots \ar 0 \ar X_n \ar X_{n-1} \ar X_{n-2} \ar
\cdots \ar X_{-n} \ar B_{-n-1} \ar 0 \ar \cdots$$ Now
Fact~\ref{kappa-generated-characterization} tells us $X = S_n$ for
some $n$. So $X$ is bounded.

For the converse assume that each $X_n$ is $\kappa$-generated
where $\kappa > \omega$. Let $X = \sum_{i \in I} S_i$ where each
$S_i$ is a subcomplex of $X$. Then $X_n = \sum_{i \in I} (S_i)_n$
and so there exists a set $J_n \subseteq I$ such that $|J_n| <
\kappa$ and $X_n = \sum_{i \in J_n} (S_i)_n$. Let $J = \bigcup_{n
\in \Z} J_n$. Then $X = \sum_{i \in J} S_i$ and $|J| < \kappa$ as
desired. For the case of when $\kappa = \omega$ it is clear that
we need the complex $X$ to be bounded so that the collection $J$
is finite.
\end{proof}

\begin{proposition}\label{prop-F-twiddle locally cogenerated}
Let $\class{F}$ be a $\kappa$-Kaplansky class. If $\cat{G}$ is
locally $\kappa$-generated, then the class $\lclass{F}$ in
$\textnormal{Ch}(\cat{G})$ is locally $\kappa$-cogenerated, as
long as $\kappa > \omega$. (If $\kappa = \omega$, then
$\lclass{F}$ is still $\kappa'$-cogenerated for any regular
$\kappa' > \omega$.)
\end{proposition}

\begin{proof}
$\lclass{F}$ is the class of all exact chain complexes $F$ with
$Z_n F \in \class{F}$ for all $n$. Suppose $0 \neq F \in
\class{F}$ is given. We wish to construct a nonzero exact complex
$S \subseteq F$ in such a way that each $Z_n S , Z_n F/Z_n S \in
\class{F}$ and each $S_n$ is $\kappa$-generated. It follows that
$\lclass{F}$ is locally $\kappa$-cogenerated if $\kappa > \omega$.
If $\kappa = \omega$ then $S$ may not be locally
$\omega$-presentable because it need not be bounded. However in
this case we still have $S$ $\kappa'$-generated for any regular
$\kappa' > \kappa$. So $\lclass{F}$ is locally
$\kappa'$-cogenerated.

Since $F \neq 0$ there must be some integer $n$ for which $Z_n F
\neq 0$. We start by finding $0 \neq S'_n \subseteq Z_n F$ with
$S'_n$ $\kappa$-presentable and $S'_n \in \class{F} , \ Z_n F/S'_n
\in \class{F}$. We set $S_i = 0$ for all $i < n$ and set $S_n =
S'_n$. We now want to inductively build $S_i$ for $i > n$.

Since $\cat{G}$ is locally $\kappa$-generated, we can use
Lemma~\ref{surjective-lemma} to find a $\kappa$-generated
subobject $X_{n+1} \subseteq F_{n+1}$ for which $d_{|X_{n+1}}
\mathcolon X_{n+1} \ar S'_n$ is an epimorphism. Since $S'_n$ is
$\kappa$-presentable and $X_{n+1}$ is $\kappa$-generated and
$\cat{G}$ is locally $\kappa$-generated, it follows from
Lemma~\ref{kappa-prentable-characterization} that
$\ker{d_{|X_{n+1}}}$ is $\kappa$-generated as well. Now using the
definition of Kaplansky class, we can find a $\kappa$-presentable
$S'_{n+1}$ such that $\ker{d_{|X_{n+1}}} \subseteq S'_{n+1}
\subseteq Z_{n+1} F$ and $S'_{n+1} , \ Z_{n+1} F / S'_{n+1} \in
\class{F}$. Now
$$d_{|X_{n+1} + S'_{n+1}} \mathcolon X_{n+1} + S'_{n+1} \ar S'_n$$ is
an epimorphism whose kernel is $S'_{n+1}$. So we set $S_{n+1} =
X_{n+1} + S'_{n+1}$. In this way we continue inductively to
construct an exact subcomplex $0 \neq S \subseteq F$ with $Z_n S ,
Z_n F / Z_n S \in \class{F}$ and $Z_n S$ $\kappa$-presentable (and
therefore $S_n$ is $\kappa$-generated by Facts~\ref{fact-b}
and~\ref{fact-b.5}).
\end{proof}

\begin{theorem}\label{them-model structure}
Let $\cat{G}$ be a locally $\kappa$-presentable Grothendieck
category. Suppose $\class{F}$ is a class of objects which
satisfies the following:

(1) $\class{F}$ is a $\kappa$-Kaplansky class.

(2) $\class{F}$ contains a $\kappa$-presentable generator $G$ for
$\cat{G}$.

(3) $\class{F}$ is closed under transfinite extensions and
retracts.

(4) $\ldgclass{F} \cap \class{E} = \lclass{F}$.

Then we have an induced model category structure on
$\textnormal{Ch}(\cat{G})$ where the weak equivalences are the
homology isomorphisms. The cofibrations (resp. trivial
cofibrations) are the monomorphisms whose cokernels are in
$\ldgclass{F}$ (resp. $\lclass{F})$. The fibrations (resp. trivial
fibrations) are the epimorphisms whose kernels are in
$\rdgclass{C}$ (resp. $\rclass{C}$), where $\class{C} =
\rightperp{\class{F}}$. Furthermore this model structure is
cofibrantly generated. The generating cofibrations form the set
$$I = \{\, 0 \ar D^n(G) \, \} \cup \{\, S^{n-1}(G) \ar D^n(G) \,\}
\cup \{\, S^n(A) \xrightarrow{} S^n(B) \, \}$$ where $A \ar B$
ranges over all possible monomorphisms with $B$ a
$\kappa$-generated object for which $B/A \in \class{F}$. The
generating trivial cofibrations are $$J = \{\, 0 \ar D^n(G) \, \}
\cup \{\, X \ar Y \,\}$$ where $X \ar Y$ ranges over all
monomorphisms where $Y$ is a $\kappa$-generated complex and $Y/X
\in \lclass{F}$
\end{theorem}

\begin{proof}
All of the work has been done. Since $\class{F}$ is a locally
cogenerated class, closed under transfinite extensions and
retracts, and containing a generator for $\cat{G}$,
Proposition~\ref{prop-locally cogeneratd class induces small
cotorsion pair} tells us $(\class{F},\class{C})$ is a small
cotorsion pair. Now by Corollary~3.8 of~\cite{gillespie} we have
the induced cotorsion pairs $(\lclass{F},\rdgclass{C})$ and
$(\ldgclass{F},\rclass{C})$ of chain complexes. These cotorsion
pairs are compatible by Corollary~\ref{cor-equivalence of
compatible and coresolving} above.

Having the cotorsion pairs of complexes $(\ldgclass{F} \cap
\class{E}, \rdgclass{C})$ and $(\ldgclass{F}, \rdgclass{C} \cap
\class{E})$ at hand we are in position to use Hovey's Theorem~2.2
from~\cite{hovey}. By Proposition~\ref{prop-induced cotorsion pair
of complexes is small} we know that $(\ldgclass{F}, \rdgclass{C}
\cap \class{E})$ is small, so this cotorsion pair is complete by
Corollary~6.6 of~\cite{hovey}. Next $\ldgclass{F} \cap \class{E}$
is locally cogenerated by Proposition~\ref{prop-F-twiddle locally
cogenerated}, so $(\ldgclass{F} \cap \class{E},\rclass{C})$ is
small too and hence complete. Theorem~2.2 and Lemma~6.7
of~\cite{hovey} tell us we have an induced model structure on
$\textnormal{Ch}(\cat{G})$ which is cofibrantly generated. One can
see that the weak equivalences, fibrations and cofibrations are as
stated in the theorem by looking at Section~5 of~\cite{hovey}.
\end{proof}

\section{Monoidal model structures}\label{sec-monoidal-model
structures}

Suppose our ground category $\cat{G}$ has a tensor product making
it a closed symmetric monoidal Grothendieck category. We would
like to have a condition on our Kaplansky class $\class{F}$ which
will guarantee that the model structure induced from
Theorem~\ref{them-model structure} will be monoidal. We refer the
reader to Chapter~4 of~\cite{hovey-model-categories} for a
detailed discussion of monoidal model structures. There is also a
nice discussion of the essential ideas in Section~7
of~\cite{hovey}. In fact our proof below will rely on basic
results from Section~7 of~\cite{hovey}. Also, the crucial argument
of the proof below, part (ii), is due to Hovey.

\begin{theorem}\label{them-monoidal condition}
Suppose $\cat{G}$ is a Grothendieck category with a closed
symmetric monoidal structure $- \tensor_{\cat{G}} -$ and let
$\class{F}$ be a class of objects such that $\cat{G}$ and
$\class{F}$ satisfy the hypotheses of Theorem~\ref{them-model
structure}. Then the induced model structure on
$\textnormal{Ch}(\cat{G})$ is monoidal with respect to the usual
tensor product of chain complexes if the following conditions
hold:

(1) Each object in $\class{F}$ is flat. I.e., $F \tensor_{\cat{G}}
-$ is exact for each $F \in \class{F}$.

(2) If $F$ and $G$ both belong to $\class{F}$, then $F
\tensor_{\cat{G}} G$ also belongs to $\class{F}$.

(3) $U \in \class{F}$, where $U$ is the unit for the monoidal
structure on $\cat{G}$.
\end{theorem}

\begin{proof}
Suppose the three stated conditions hold. We wish to check the
four conditions of Theorem~7.2 from~\cite{hovey}. In this case the
four conditions translate to the following: \vspace{.1in}

(i) Every cofibration is a pure injection in each degree.

(ii) If $X$ and $Y$ are in $\ldgclass{F}$, then $ X \tensor Y$ is
  in $\ldgclass{F}$.

(iii) If $X$ is in $\lclass{F}$ and $Y$ is in $\ldgclass{F}$, then
$X \tensor Y$ is
  in $\lclass{F}$.

(iv) The unit for the monoidal structure on Ch($\cat{G}$) is in
$\ldgclass{F}$. \vspace{.1in}

Before we check each condition we make a few observations. First,
as we can see from Proposition~\ref{prop-induced cotorsion pair of
complexes is small}, the hypotheses guarantee that
$(\ldgclass{F},\rclass{C})$ is a small cotorsion pair with
cogenerating set consisting of spheres on objects from
$\class{F}$. That is, the cogenerating set consists of complexes
of the form $S^n(F)$ where $F \in \class{F}$. For the rest of the
proof we will call any such object $S^n(F)$ an
``$\class{F}$-sphere''. Second, $\ldgclass{F}$ contains a set of
generators for Ch($\cat{G}$) since if $G \in \class{F}$ generates
$\cat{G}$, then the set of complexes $D^n(G)$ generates
Ch($\cat{G}$). (Each $D^n(F)$ is in $\ldgclass{F}$ by Lemma~3.4
of~\cite{gillespie}.) It now follows from
Lemma~\ref{lemma-classification of left side of cotorsion pair}
that each complex in $\ldgclass{F}$ is a direct summand of a
transfinite extension of $\class{F}$-spheres. Lastly, we note that
if $X \in \ldgclass{F}$, then $X \tensor -$ is exact. Indeed one
can easily check that a complex $X$ is tensor exact if and only if
$X_n$ is flat for each $n$. We now check the four conditions
above.

(i) This follows from a basic property of flat objects.
Cofibrations are monomorphisms with cokernels in $\ldgclass{F}$.
In particular, the cokernel is flat in each degree. Therefore each
cofibration is a pure injection in each degree. (See for example,
the proof of Lemma~{XVI.3.1} in~\cite{lang}.)

(ii) \emph{Step 1:} We first show that if $X$ is a transfinite
extension of $\class{F}$-spheres and if $S^n(F)$ is an
$\class{F}$-sphere, then $X
\tensor_{\textnormal{Ch}(\cat{G})}S^n(F)$ is in $\ldgclass{F}$. To
see this first notice that the tensor product of two
$\class{F}$-spheres is again an $\class{F}$-sphere and in
particular is in $\ldgclass{F}$. Now say $X$ is a transfinite
extension of a sequence such as $$X_0 \hookrightarrow X_1
\hookrightarrow \cdots \hookrightarrow X_{\alpha} \hookrightarrow
X_{\alpha + 1} \hookrightarrow \cdots$$ where $X_0$ and each
$X_{\alpha + 1}/X_{\alpha}$ are $\class{F}$-spheres. Since the
functor $- \tensor_{\textnormal{Ch}(\cat{G})}S^n(F)$ is exact and
also preserves direct limits, applying it to the above sequence
will display $X \tensor_{\textnormal{Ch}(\cat{G})}S^n(F)$ as a
transfinite extension of objects in $\ldgclass{F}$. Therefore it
too is in $\ldgclass{F}$ by Lemma~6.2 of~\cite{hovey}.

\emph{Step 2:} Now we see that if $X \in \ldgclass{F}$ and if
$S^n(F)$ is an $\class{F}$-sphere, then $X
\tensor_{\textnormal{Ch}(\cat{G})}S^n(F)$ is in $\ldgclass{F}$.
The reason is that such an $X$ is a direct summand of a
transfinite extension of $\class{F}$-spheres. Since the functor $-
\tensor_{\textnormal{Ch}(\cat{G})}S^n(F)$ commutes with direct
sums, the result follows from the fact that $\ldgclass{F}$ is
closed under direct summands. (The left side of any cotorsion pair
is closed under direct summands.)

\emph{Step 3:} Now we argue the same way to see that $X
\tensor_{\textnormal{Ch}(\cat{G})} Y$ is in $\ldgclass{F}$
whenever $X$ and $Y$ are in $\ldgclass{F}$. For example if $Y$ is
in $\ldgclass{F}$ then it is a retract of a transfinite extension
of $\class{F}$-spheres. As above, since $X
\tensor_{\textnormal{Ch}(\cat{G})} -$ commutes with direct sums
and $\ldgclass{F}$ is closed under taking direct summands we just
need to show that $X \tensor_{\textnormal{Ch}(\cat{G})} Y$ is in
$\ldgclass{F}$ when $X$ is in $\ldgclass{F}$ and $Y$ is a
transfinite extension of $\class{F}$-spheres. But in this case $Y$
can be viewed as the transfinite extension of a sequence such as
$$Y_0 \hookrightarrow Y_1 \hookrightarrow \cdots \hookrightarrow
Y_{\alpha} \hookrightarrow Y_{\alpha + 1} \hookrightarrow \cdots$$
where $Y_0$ and each $Y_{\alpha + 1}/Y_{\alpha}$ are
$\class{F}$-spheres. Then we apply the functor $X
\tensor_{\textnormal{Ch}(\cat{G})} -$ and use the fact from Step 2
to argue that $X \tensor_{\textnormal{Ch}(\cat{G})} Y$ is a
transfinite extension of complexes in $\ldgclass{F}$. So $X
\tensor_{\textnormal{Ch}(\cat{G})} Y$ is in $\ldgclass{F}$ too.

(iii) The approach is similar to our proof of (ii). Note that if
$\class{E}$ is the class of exact complexes in Ch($\cat{G}$), then
$\class{E}$ is closed under transfinite extensions and direct
summands. One reason this is true is that $\class{E}$ is the left
side of the cotorsion pair $(\class{E}, \rdgclass{I})$ where
$\rdgclass{I}$ is the class of dg-injective complexes.

\emph{Step 1:} First we will prove that if $E$ is an exact complex
in Ch($\cat{G}$) and $Y \in \ldgclass{F}$, then $E \tensor Y$ is
exact. We leave it to the reader to see that, as in our proof of
(ii), we can assume that $Y$ is a transfinite extension of
$\class{F}$-spheres. So suppose $Y$ is the direct limit of a
sequence such as
$$Y_0 \hookrightarrow Y_1 \hookrightarrow \cdots \hookrightarrow
Y_{\alpha} \hookrightarrow Y_{\alpha + 1} \hookrightarrow \cdots$$
where $Y_0$ and each $Y_{\alpha + 1}/Y_{\alpha}$ are
$\class{F}$-spheres. Since each short exact sequence $$0 \ar
Y_{\alpha} \ar Y_{\alpha + 1} \ar Y_{\alpha + 1}/Y_{\alpha} \ar
0$$ is pure in each degree, the sequence below is also exact:
$$0 \ar E \tensor Y_{\alpha} \ar E \tensor Y_{\alpha + 1} \ar E
\tensor Y_{\alpha + 1}/Y_{\alpha} \ar 0$$ Now by applying the
functor $E \tensor -$ to the entire sequence
$$Y_0 \hookrightarrow Y_1 \hookrightarrow \cdots \hookrightarrow
Y_{\alpha} \hookrightarrow Y_{\alpha + 1} \hookrightarrow \cdots$$
we can argue that $E \tensor Y$ is a transfinite extension of
exact complexes and so is itself exact.

\emph{Step 2:} Now say $X \in \lclass{F}$ and $Y \in
\ldgclass{F}$. Then $X \in \ldgclass{F}$ by our assumption that
$\lclass{F} = \ldgclass{F} \cap \class{E}$. So by part (ii) we
have that $X \tensor Y$ is in $\ldgclass{F}$. But by Step 1 of
this proof we also have that $X \tensor Y$ is in $\class{E}$.
Therefore $X \tensor Y$ is in $\lclass{F}$.

(iv) Since $U \in \class{F}$, we have from Lemma~3.4
of~\cite{gillespie} that $S(U) = S^0(U)$ is in the class
$\ldgclass{F}$. Since $S(U)$ is the the unit for the tensor
product on Ch($\cat{G}$) we are done.
\end{proof}


\section{Chain complexes of quasi-coherent sheaves}\label{sec-quasi-coherent sheaves}

For the rest of the paper we turn to applications of
Theorem~\ref{them-model structure} and Theorem~\ref{them-monoidal
condition}. In this section we prove the existence of the flat
model structure on Ch($\qcox$) when $X$ is a quasi-compact and
semi-separated scheme.

So let $(X, \class{O}_X)$ be a scheme on a topological space $X$.
We denote the category of all sheaves of $\class{O}_X$-modules by
$\class{O}_X$-Mod. We denote by $\qcox$, the full subcategory of
$\class{O}_X$-Mod consisting of all quasi-coherent
$\class{O}_X$-modules. The class $\class{F}$ of all flat
quasi-coherent sheaves will be the Kaplansky class inducing the
model structure. As far as the author can tell we need some
assumptions on the scheme $X$ in order to satisfy the hypotheses
of Theorem~\ref{them-model structure}. From a study of the current
literature it seems as though assuming $X$ is quasi-compact and
semi-separated is the best we can do. In this case we see
from~\cite{lipman} or~\cite{murfet} that every quasi-coherent
$\class{O}_X$-module is the quotient of a flat quasi-coherent
$\class{O}_X$-modules. We start this section by recalling some
facts about the category $\qcox$ as well as defining
semi-separated. The author has used~\cite{hartshorne},
\cite{litaka} and Appendix~B of~\cite{thomason-trobaugh} as basic
references on this material. We will also use recent results
from~\cite{enochs-quasi-coherent}.

The inclusion functor $\qcox \ar \class{O}_X$-Mod is exact and
$\qcox$ is an abelian subcategory of $\class{O}_X$-Mod. Therefore
finite limits and colimits in $\qcox$ are taken as in
$\class{O}_X$-Mod. As stated in the beginning of Appendix~B
of~\cite{thomason-trobaugh}, $\qcox$ is cocomplete with all small
colimits taken as in $\class{O}_X$-Mod. Therefore direct limits
are exact in $\qcox$. In fact $\qcox$ also has a set of generators
making it a Grothendieck category. This last fact was stated as
Lemma~2.1.7 in~\cite{conrad} and apparently is due to Gabber. A
proof is not given in~\cite{conrad} but an independent proof can
be found in~\cite{enochs-quasi-coherent}.

Since we know $\qcox$ has a generator and the inclusion functor
$\qcox \ar \class{O}_X$-Mod preserves small colimits the inclusion
has a right adjoint $Q \mathcolon \class{O}_X$-Mod $\ar \qcox$ by
the special adjoint functor theorem. Thus $\qcox$ is a
coreflective subcategory of $\class{O}_X$-Mod. In particular, $Q$
satisfies a universal property dual in nature to the
sheafification of a presheaf. One can now easily check that all
small limits exist in $\qcox$ and are taken by applying $Q$ to the
usual limit taken in $\class{O}_X$-Mod. Lastly, we recall that
$\qcox$ is closed under extensions and tensor products taken in
$\class{O}_X$-Mod. The tensor product makes $\qcox$ into a closed
symmetric monoidal category. The closed structure is given by
applying the functor $Q$ after the usual ``sheafhom'' functor.

A quasi-coherent $\class{O}_X$-module $F$ is flat if the tensor
functor $F \tensor_{\class{O}_X} -$ is exact. Clearly the class of
flat modules in $\qcox$ is merely the class of flat modules in
$\class{O}_X$-Mod intersected with the class of all quasi-coherent
$\class{O}_X$-modules. We denote the class of all flat
quasi-coherent $\class{O}_X$-modules by $\class{F}$. We point out
that since flatness is defined in terms of the tensor product
which is a left adjoint, $\class{F}$ is closed under direct limits
and retracts (summands). Furthermore, if $$0 \ar F' \ar F \ar F''
\ar 0$$ is a short exact sequence of quasi-coherent sheaves with
$F'' \in \class{F}$, then $F' \in \class{F}$ if and only if $F \in
\class{F}$. These assertions are all true by using categorical
arguments with the tensor product, similar to those found
in~\cite{lang} for modules over a ring. Alternatively, they can be
shown using the characterization of flatness as a ``stalkwise''
property along with the fact that the corresponding property holds
for modules over a ring. In any case, we point out that
$\class{F}$ is closed under retracts and transfinite extensions
and also that $\class{F}$ is resolving, meaning it is closed under
taking kernels of epimorphisms between objects in $\class{F}$.

We say a scheme $X$ is semi-separated if there exists an affine
basis $\class{V} = \{V_{\alpha}\}$ for sp($X$) which is closed
under finite intersections. The basis $\class{V}$ is called a
\emph{semi-separating affine basis} of $X$. As pointed out in~B.7
of~\cite{thomason-trobaugh} a semi-separated scheme is
quasi-separated. Furthermore, a separated scheme is semi-separated
with the set of all affine subsets serving as a semi-separating
affine basis.

It is well-known that both $\cat{O}_X$-Mod and $\qcox$ do not, in
general, have enough projectives. While the category
$\cat{O}_X$-Mod has a set of flat generators, the author does not
know whether or not for a general scheme $X$, $\qcox$ has a set of
flat generators, or equivalently whether each quasi-coherent sheaf
can be written as the quotient of a flat quasi-coherent sheaf.
However, Proposition~1.1 of~\cite{lipman} implies that this is
indeed the case if we assume $X$ is a quasi-compact, separated
scheme. In fact as pointed out by Daniel Murfet their proof works
for any quasi-compact semi-separated scheme $X$. Murfet's proof
can be found as Proposition~16 in~\cite{murfet}. Since the
quasi-compact semi-separated hypothesis is the best one the author
has found for the existence of enough flats in $\qcox$ we will
build the flat model structure on Ch($\qcox$) under these
assumptions on $X$.

The following subsections will contain lemmas breaking down the
proof of Theorem~\ref{them-main} which of course amounts to
checking the four hypotheses of Theorem~\ref{them-model structure}
for the class $\class{F}$ of flat modules in $\qcox$. We start
with the Kaplansky class condition. As we will see, it readily
follows from recent work of Enochs and Estrada that $\class{F}$ is
indeed a Kaplansky class for \emph{any} scheme
$X$~\cite{enochs-quasi-coherent}. Since the term \emph{Kaplansky
class} was not used in~\cite{enochs-quasi-coherent} we now check
carefully that their work proves $\class{F}$ is a Kaplansky class
in the sense of Definition~\ref{def-kaplansky class}.


\subsection{Quasi-coherent modules on a representation of a quiver.}\label{subsec-quivers}

By a \emph{quiver} we mean a directed graph $Q = (V,E)$ where $V$
is the set of vertices and $E$ is the set of edges. By a
\emph{representation} of $Q$, we will mean a functor $R \mathcolon
Q \ar \textnormal{CRng}$ where $Q$ is thought of as a category in
the obvious way and $\textnormal{CRng}$ is the category of
commutative rings with identity. Let $R$ be a representation of
$Q$. By an $R$-module $M$ we will mean an $R(v)$-module, $M(v)$,
for each $v \in V$, and an $R(v)$-linear map $M(e) \mathcolon M(v)
\ar M(w)$ for each edge $e \mathcolon v \ar w$ in $E$. We refer
the reader to Section~2 of~\cite{enochs-quasi-coherent} for the
definition of a \emph{flat representation} and the definition of a
\emph{quasi-coherent R-module}. It is also pointed out there that
the category $\cat{C}$ of quasi-coherent $R$-modules is a
Grothendieck category when $R$ is flat.

We define the cardinality of an $R$-module as $$|M| = |\coprod_{v
\in V} M(v)|.$$

\begin{lemma}\label{presentable-generated-cardinality-equal}
Let $R$ be a flat representation of a quiver $Q = (V,E)$ and $M$
be a quasi-coherent $R$-module. Let $\kappa$ be a regular cardinal
for which $\kappa > |R(v)|$ for all $v$ and $\kappa > \text{max}
\{|E|,|V|\}$. Then the following are equivalent.

(1) $|M| < \kappa$.

(2) $M$ is $\kappa$-generated.

(3) $M$ is $\kappa$-presentable.
\end{lemma}

\begin{proof}
(1) $\Rightarrow (2)$. We use
Fact~\ref{kappa-generated-characterization}. Suppose $|M| <
\kappa$ and say $M = \sum_{i \in I} M_i$ is a $\kappa$-filtered
union of quasi-coherent submodules. For each $x \in M$, there
corresponds some $i \in I$ such that $x \in M_i$. But $\sum_{i \in
I} M_i$ is $\kappa$-filtered and $|M| < \kappa$, so there exists
$i \in I$ such that $M = M_i$.

$(2) \Rightarrow (1)$. Let $\class{S}$ be the collection of
\emph{all} subsets $S \subseteq \coprod_{v \in V} M(v)$ such that
$|S| < \kappa$. For each $S\in \class{S}$, let $M_S$ represent the
quasi-coherent submodule generated by $S$. Then $|M_S| < \kappa$.
(One way to see this is to use Proposition~3.3
of~\cite{enochs-quasi-coherent}. All $\leq$ signs in that
Proposition can be changed to $<$ signs.) Note that $(\class{S},
\subseteq)$ is $\kappa$-filtered and in fact $M$ is the
$\kappa$-filtered union $\sum_{S \in \class{S}} M_S$. By
Fact~\ref{kappa-generated-characterization}, $M = M_S$ for some $S
\in \class{S}$.

$(3) \Rightarrow (2)$ is automatic. We now prove $(2) \Rightarrow
(3)$, using that (1) iff (2). First we point out that the category
of quasi-coherent $R$-modules is locally $\kappa$-generated.
Indeed each $M$ can be expressed as the $\kappa$-filtered union $M
= \sum_{S \in \class{S}} M_S$ where each $M_S$ is
$\kappa$-generated as in the last paragraph. Therefore, we may use
the characterization of $\kappa$-presentable objects in
Fact~\ref{kappa-prentable-characterization}. Suppose $M$ is
$\kappa$-generated and $N \ar M$ is an epimorphism with $N$ a
$\kappa$-generated $R$-module. Then $|N| < \kappa$. So of course
$|\ker{(N \ar M)}| < \kappa$ (kernels are computed componentwise),
which means $\ker{(N \ar M)}$ is $\kappa$-generated. This proves
$M$ is $\kappa$-presentable.
\end{proof}

An $R$-module $M$ is called \emph{flat} if each $M(v)$ is a flat
$R(v)$-module.

\begin{lemma}\label{lemma-flat quasi-coherents are kaplansky}
Let $R$ be a flat representation of a quiver $Q = (V,E)$. Let
$\kappa$ be a regular cardinal for which $\kappa
> |R(v)|$ for all $v$ and $\kappa > \text{max} \{|E|,|V|\}$. Then
the class of all flat quasi-coherent $R$-modules constitute a
$\kappa$-Kaplansky class.
\end{lemma}

\begin{proof}
This follows immediately from
Lemma~\ref{presentable-generated-cardinality-equal} along with the
Proposition~3.3 of~\cite{enochs-quasi-coherent}. Again, all $\leq$
signs can be changed to $<$ signs in the statement and proof of
Proposition~3.3 of~\cite{enochs-quasi-coherent}.
\end{proof}


\subsection{A category equivalent to
$\textbf{Qco(X)}$}\label{subsec-category equivalent to Qco(X)}

If $X$ is a scheme, then the collection of all affine open subsets
determines a directed graph, $\class{A}_X$. The opposite graph
(reverse all arrows) is also a quiver which we will denote $Q_X$.
We get a (flat) representation $R$ by letting $R(U) =
\class{O}_X(U)$ for each $U \in \class{A}_X$ and using the
restriction maps to get $R(V) \ar R(U)$ whenever $U \subseteq V$.
In the same way, a quasi-coherent sheaf $S$ on $X$ gives rise to a
quasi-coherent $R$-module $M$ by letting $M(U) = S(U)$ for $U \in
\class{A}_X$. In fact, as explained
in~\cite{enochs-quasi-coherent}, the category $\cat{C}$ of
quasi-coherent $R$-modules is categorically equivalent to $\qcox$.
Furthermore, the equivalence restricts to an equivalence between
the flat objects in $\cat{C}$ and the flat objects of $\qcox$.

\begin{lemma}\label{lemma-equivalence preserves a kaplansky class}
Let $\class{A}_0 \subseteq \cat{A}$ and $\class{B}_0 \subseteq
\cat{B}$ where $\cat{A}$ and $\cat{B}$ are abelian categories. Let
$G \mathcolon \cat{A} \ar \cat{B}$ be an additive functor which is
an equivalence. If $G$ restricts to an equivalence between
$\class{A}_0$ and $\class{B}_0$, then $\class{A}_0$ is a
$\kappa$-Kaplansky class if and only if $\class{B}_0$ is a
$\kappa$-Kaplansky class.
\end{lemma}

\begin{proof}
Since $G \mathcolon \cat{A} \ar \cat{B}$ is an equivalence, there
exists $H \mathcolon \cat{B} \ar \cat{A}$ such that $HG \cong
1_{\cat{A}}$ and $GH \cong 1_{\cat{B}}$. Furthermore, $G$ is both
a left and right adjoint of $H$ (and so each preserves colimits).
We start by showing that $G$ preserves $\kappa$-generated and
$\kappa$-presentable objects. Indeed if $X \in \cat{A}$ is
$\kappa$-presentable then the functor $\cat{A}(X,H(-))$ preserves
$\kappa$-filtered colimits. But $\cat{A}(X,H(-)) \cong
\cat{B}(G(X),-)$, so $G(X)$ must be $\kappa$-presentable. On the
other hand, if $G(X)$ is $\kappa$-presentable, then
$\cat{B}(G(X),G(-))$ preserves $\kappa$-filtered colimits.
Therefore, $\cat{B}(G(X),G(-)) \cong \cat{A}(HG(X),-)$ preserves
$\kappa$-filtered colimits and $HG(X) \cong X$ must be
$\kappa$-presentable. The same proof (but with $\kappa$-filtered
colimits of monomorphisms) shows that an equivalence preserves
$\kappa$-generated objects.

Now we suppose $G$ restricts to an equivalence between
$\class{A}_0$ and $\class{B}_0$ and that $\class{B}_0$ is a
$\kappa$-Kaplansky class of $\cat{B}$. We will show that
$\class{A}_0$ is a $\kappa$-Kaplansky class of $\cat{A}$. So let
$X \subseteq F \neq 0$ where $F \in \class{A}_0$ and $X$ is a
$\kappa$-generated object in $\cat{A}$. Then $G(X) \subseteq G(F)
\neq 0$ and $G(X)$ is $\kappa$-generated and $G(F) \in
\class{B}_0$. Since $\class{B}_0$ is a $\kappa$-Kaplansky class,
there exists a nonzero $\kappa$-presentable object $S \in
\class{B}_0$ for which $G(X) \subseteq S \subseteq G(F)$ and
$G(F)/S \in \class{B}_0$. Therefore $H(S) \in \class{A}_0$ is a
nonzero $\kappa$-presentable object for which $X = HG(X) \subseteq
H(S) \subseteq HG(F) = F$ and $F/H(S) = HG(F)/H(S) = H(GF/S) \in
\class{A}_0$.
\end{proof}

\begin{proposition}\label{prop-flat cotorsion pair}
Let $(X , \class{O}_X)$ be quasi-compact and semi-separated
scheme. Let $\class{F}$ be the class of flat quasi-coherent
$\class{O}_X$-modules and define $\class{C} =
\rightperp{\class{F}}$. Then $\class{F}$ is a Kaplansky class and
$(\class{F} , \class{C})$ is a small cotorsion pair.
\end{proposition}

\begin{proof}
It follows from the lemmas above that $\class{F}$ is a Kaplansky
class. In particular $\class{F}$ is locally cogenerated.
Furthermore as explained in the introduction to this section,
$\class{F}$ is closed under transfinite extensions and retracts
and contains a generator. Therefore $(\class{F} , \class{C})$ is a
small cotorsion pair by Proposition~\ref{prop-locally cogeneratd
class induces small cotorsion pair}.
\end{proof}


\subsection{Exact dg-flat complexes are
flat}\label{subsec-compatibility of flats}

Let $(X , \class{O}_X)$ be a quasi-compact semi-separated scheme.
For what follows we fix a semi-separating affine basis $\class{V}$
and again let $\class{F}$ be the class of all flat modules in
$\qcox$. We will prove that $\lclass{F} = \ldgclass{F} \cap
\class{E}$ where $\class{E}$ is the class of all exact chain
complexes.

Recall that if $f \mathcolon Y \ar Z$ is any morphism of schemes
then by Proposition~II.5.8 of~\cite{hartshorne} the inverse image
functor $f^*$ preserves quasi-coherence. In particular, if
$V_{\alpha} \in \class{V}$ and $j \mathcolon V_{\alpha} \ar X$ is
the inclusion, then we have the functor $j^* \mathcolon
\textnormal{Qco}(X) \ar \textnormal{Qco}(X|_{V_{\alpha}})$. In
this case the functor is merely restriction of a quasi-coherent
$\class{O}_X$-module to the affine subset $V_{\alpha}$ and it is
clear that $j^*$ is an exact functor.

On the other hand, recall that the direct image functor $f_*$ does
not always preserve quasi-coherence. However it is easy to see
that $j_*$ does preserve quasi-coherence. Indeed $j$ is separated
since it is an open immersion (see Theorem~1.17 of~\cite{litaka}).
Also since affine subsets are quasi-compact we see by
Proposition~1.51 of~\cite{litaka} (with our semi-separating basis
$\class{V}$ used as the affine open cover required in the
Proposition), that $j_*$ preserves quasi-coherence.

\begin{lemma}\label{lemma-j_* is exact and preserves cotorsions}
Let $(X, \class{O}_X)$ be a semi-separated scheme with
semi-separating affine basis $\class{V} = \{V_{\alpha}\}$. Let
$V_{\alpha} \in \class{V}$ and $j \mathcolon V_{\alpha} \ar X$ be
the inclusion. Then $j^* \mathcolon \textnormal{Qco}(X) \ar
\textnormal{Qco}(X|_{V_{\alpha}})$ is left adjoint to $j_*
\mathcolon \textnormal{Qco}(X|_{V_{\alpha}}) \ar
\textnormal{Qco}(X)$. Furthermore $j_*$ is exact and preserves
cotorsion objects.
\end{lemma}

\begin{proof}
It is a standard fact that inverse image functors are left adjoint
to direct image functors. For example, see Section~II.5
of~\cite{hartshorne}. In particular, $j^*$ is left adjoint to
$j_*$. Since $j_*$ is a right adjoint it is left exact. We will
show that $j_*$ preserves surjections. So suppose we have a
surjection $\class{F} \ar \class{G}$ of quasi-coherent sheaves in
$\textnormal{Qco}(X|_{V_{\alpha}})$. Then given any affine subset
$V_{\gamma} \subseteq V_{\alpha}$, we know $\class{F}(V_{\gamma})
\ar \class{G}(V_{\gamma})$ is a surjection. It follows that for
any affine $V_{\beta} \in \class{V}$, $[j_*(\class{F})](V_{\beta})
\ar [j_*(\class{G})](V_{\beta})$ is surjective. Indeed the
definition of $[j_*(\class{F})](V_{\beta}) \ar
[j_*(\class{G})](V_{\beta})$ is just $\class{F}(V_{\gamma}) \ar
\class{G}(V_{\gamma})$ where $V_{\gamma} = V_{\alpha} \cap
V_{\beta} \in \class{V}$. It follows that $j_*$ is exact since
$\class{V}$ is a basis for the topology on $X$.

We now show that for any quasi-coherent sheaf $\class{F} \in
\textnormal{Qco}(X)$ and any quasi-coherent sheaf $\class{C} \in
\textnormal{Qco}(X|_{V_\alpha})$ we have an isomorphism
$$\Ext^n(j^* \class{F},\class{C}) \cong \Ext^n(\class{F} ,j_*
\class{C}).$$ First note that $j_*$ being the right adjoint of an
exact functor preserves injective objects by Proposition~2.3.10
of~\cite{weibel}. Since $j_*$ itself is also exact it preserves
injective resolutions. The result now follows by taking an
injective resolution of $\class{C}$, applying $j_*$ to the
resolution, and then applying the adjoint relationship between
$j_*$ and $j^*$.

Now suppose $\class{C} \in \textnormal{Qco}(X|_{V_{\alpha}})$ is
cotorsion and $\class{F}$ is a flat quasi-coherent
$\class{O}_X$-module. We wish to show that $j_* \class{C}$
cotorsion. But $j^* \class{F}$ is clearly flat since $j^*$ is just
restriction, and so $\Ext^1(\class{F}, j_* \class{C}) \cong
\Ext^1(j^* \class{F},\class{C}) = 0$. Therefore $j_* \class{C}$ is
cotorsion.
\end{proof}

\begin{proposition}\label{prop-flat quasi-coherent sheaf classes
are compatible} Let $(X, \class{O}_X)$ be a semi-separated scheme
with semi-separating affine basis $\class{V} = \{V_{\alpha}\}$ and
let $\class{F}$ be the class of flat quasi-coherent
$\class{O}_X$-modules. Then $\lclass{F} = \ldgclass{F} \cap
\class{E}$ where $\class{E}$ is the class of exact complexes.
\end{proposition}

\begin{proof}
We saw in Proposition~\ref{prop-flat cotorsion pair} that
$(\class{F},\class{C})$ is a cotorsion pair. Now Lemma~3.10
of~\cite{gillespie} says that $\lclass{F} \subseteq \ldgclass{F}
\cap \class{E}$. So we wish to prove $\lclass{F} \supseteq
\ldgclass{F} \cap \class{E}$. We let $Y$ be an exact dg-flat
complex of quasi-coherent $\class{O}_X$-modules. We want to show
that it is a flat complex. I.e. that is has flat cycles in each
degree. One way to do this is to show that the complex $j^*(Y)$ is
a flat complex of $\class{O}|_{V_{\alpha}}$-modules for each
inclusion $j \mathcolon V_{\alpha} \ar X$ where $V_{\alpha} \in
\class{V}$. So let $j \mathcolon V_{\alpha} \ar X$ be such an
inclusion. Then by the usual equivalence between the category of
$A_{\alpha}$-modules and the category of quasi-coherent sheaves on
$\tilde{A_{\alpha}}$ where $A_{\alpha}$ is a ring in which
$\tilde{A_{\alpha}} \cong \class{O}|_{V_{\alpha}}$, the chain
complex $j^*(Y) = Y|_{V_{\alpha}}$ corresponds to a chain complex
of $A_{\alpha}$-modules. Since the equivalence preserves flat
objects by Proposition~III.9.2 of~\cite{hartshorne} it also
preserves cotorsion objects. Therefore, the complex of
$A_{\alpha}$-modules corresponding to $j^*(Y)$ is flat
(respectively dg-flat) if and only if $j^*(Y)$ is flat
(respectively dg-flat). Finally, since it is already known
from~\cite{gillespie} that exact, dg-flat complexes of
$A_{\alpha}$-modules are flat complexes we will be able to
conclude that $j^*(Y)$ is flat by just proving it is exact and
dg-flat.

It is clear that $j^*(Y)$ is exact so we will show that $j^*(Y)$
is a dg-flat complex. So let $j^*(Y) \ar C$ be a morphism where
$C$ is a cotorsion complex of quasi-coherent
$\class{O}_{V_{\alpha}}$-modules. (This means that $C$ is exact
and each cycle is cotorsion.) Using the definition of dg-flat we
wish show this map is null homotopic. But using the adjoint
property of $j^*$ and $j_*$ gives us a morphism $Y \ar j_*(C)$ of
complexes. Since $j_*(-)$ is exact and preserves cotorsion modules
by Lemma~\ref{lemma-j_* is exact and preserves cotorsions} we see
that $j_*(C)$ is a cotorsion complex and by definition of dg-flat
$Y \ar j_*(C)$ is null homotopic. From this we can deduce that
$j^*(Y) \ar C$ is also null homotopic.
\end{proof}

\subsection{The main theorem}

We are now ready to prove the main theorem. By applying
Definition~3.3 of~\cite{gillespie} to the cotorsion pair
$(\class{F}, \class{C})$ we get four classes of chain complexes in
$\chqcox$. We call the complexes in $\ldgclass{F}$ the
\emph{dg-flat complexes}, and the complexes in $\lclass{F}$ the
\emph{flat complexes}. We call the complexes in
$\rdgclass{\class{C}}$ the \emph{dg-cotorsion complexes}, and the
complexes in $\rclass{\class{C}}$ the \emph{cotorsion complexes}.

\begin{theorem}\label{them-main}
Let $X$ be a quasi-compact semi-separated scheme and let
$\class{F} \subseteq \qcox$ be the class of flat quasi-coherent
sheaves. Then we have a cofibrantly generated model category
structure on $\chqcox$ which is described as follows: The weak
equivalences are the homology isomorphisms. The cofibrations
(resp. trivial cofibrations) are the monomorphisms whose cokernels
are dg-flat complexes (resp. flat complexes). The fibrations
(resp. trivial fibrations) are the epimorphisms whose kernels are
dg-cotorsion complexes (resp. cotorsion complexes). Furthermore,
this model structure is monoidal with respect to the usual tensor
product of chain complexes.
\end{theorem}

\begin{proof}
From~B.3 of~\cite{thomason-trobaugh} we see that if $X$ is
quasi-compact and quasi-separated, then $\qcox$ is locally
finitely presentable~\cite{thomason-trobaugh}. That is, $\qcox$ is
locally $\kappa$-presentable (and hence $\kappa$-generated) for
any regular cardinal $\kappa$. Since $\qcox$ is a Grothendieck
category it has a generator. By writing this generator as the
quotient of a flat quasi-coherent $\class{O}_X$-module we obtain a
flat generator. By Lemmas~\ref{lemma-flat quasi-coherents are
kaplansky} and~\ref{lemma-equivalence preserves a kaplansky class}
we can find arbitrarily large regular cardinals $\kappa$ for which
$\class{F}$ is a $\kappa$-Kaplansky class. In light of
Fact~\ref{fact-a} we can find a single regular cardinal $\kappa$
such that (i) $\qcox$ is locally $\kappa$-presentable, (ii)
$\class{F}$ is a $\kappa$-Kaplansky class, and (iii) $\class{F}$
contains a $\kappa$-presentable generator. This verifies
properties (1) and (2) of Theorem~\ref{them-model structure}.
Property (3) holds as explained in the introduction to this
section. Property (4) is Proposition~\ref{prop-flat quasi-coherent
sheaf classes are compatible}. This finishes the proof of the
existence of the flat model structure on $\chqcox$.

We now check that the model structure is monoidal. Condition (1)
of Theorem~\ref{them-monoidal condition} is trivial. Condition (2)
is true since the tensor product of two flat quasi-coherent
sheaves is again a flat quasi-coherent sheaf. The unit for the
tensor product is the structure sheaf $\class{O}_X$. It is both
quasi-coherent and flat and so condition (3) holds.
\end{proof}

Many questions come to mind that the author has not been able to
figure out or has not had time to consider. First, is there an
easier proof that $\qcox$ has enough flat objects when $X$ is
quasi-compact and semi-separated? Can we construct some actual
flat generators for $\qcox$ that are strongly dualizable in
$\cat{D}(\qcox)$? Is $\cat{D}(\qcox)$ a stable homotopy category
in the sense of~\cite{hovey-axiomatic stable homotopy}? What are
the implications of May's additivity theorem in~\cite{may3}? Can
one easily treat the derived inverse and direct image functors
from the model category viewpoint?


\section{Other applications}\label{sec-other applications}
In this section we point out other applications of
Theorem~\ref{them-model structure}. In particular we will show
that all model structures the author has previously constructed
in~\cite{gillespie} and~\cite{gillespie-sheaves} as well as the
more common ``injective'' and ``projective'' model structures on
$\ch$ can each be construed as corollaries to
Theorem~\ref{them-model structure}. Note that condition (1) (the
Kaplansky class condition) in Theorem~\ref{them-model structure}
is defined in terms of $\kappa$-generated and $\kappa$-presentable
objects. In practice it may be cumbersome to work with the
definitions of $\kappa$-generated or $\kappa$-presentable. In each
corollary below we find it convenient to pick $\kappa$ large
enough so that the $\kappa$-generated objects and
$\kappa$-presentable objects coincide. In the case when $\cat{G}$
is a concrete category we go further and relate the notions of
$\kappa$-generated and $\kappa$-presentable to the cardinality of
the underlying set. As mentioned after the definition of Kaplansky
class in Section~\ref{subsec-kaplansky classes},
$\kappa$-generated and $\kappa$-presentable become good
categorical replacements for the notion of cardinality.

\begin{corollary}
Let $\cat{G}$ be any Grothendieck category and let $\class{A}$ be
the class of all objects in $\cat{G}$. Then $\class{A}$ is a
Kaplansky class and in induces the injective model structure on
$\textnormal{Ch}(\cat{G})$. Furthermore this model structure is
cofibrantly generated.
\end{corollary}

\begin{proof}
$\cat{G}$ is locally $\kappa$-presentable for some regular
cardinal $\kappa$ by Proposition~3.10 of~\cite{beke}. Note
$\rightperp{\class{A}} = \class{I}$ is the class of injectives and
$(\class{A},\class{I})$ is the ``injective cotorsion pair''. In
order to use Theorem~\ref{them-model structure} we would like to
say that $\class{A}$ is a $\kappa$-Kaplansky class, but the author
suspects that this is just not true. (If it were, then clearly
every locally $\kappa$-generated Grothendieck category would be
locally $\kappa$-presentable. See the paragraph before
Fact~\ref{fact-c} in
Appendix~\ref{appendix-locally-presentable-cats}.) In any case
there is a trick: Pick $\kappa'$ to be a regular cardinal large
enough so that (i) $\cat{G}$ is locally $\kappa'$-presentable,
(ii) the class of $\kappa'$-generated objects coincides with the
class of $\kappa'$-presentable objects, and (iii) there is a
$\kappa'$-presentable generator for $\cat{G}$. Condition (i) is
possible by Fact\ref{fact-c.5}, condition (ii) is possible by
Appendix~\ref{appendix-kappa-generated and kappa-presentable
objects coincide}, and condition (iii) is possible by
Fact~\ref{fact-a}. Now using $\kappa'$ and $\class{A}$ we can
easily check the conditions of Theorem~\ref{them-model structure}.
Conditions (1) and (3) hold trivially and condition (2) holds by
our choice of $\kappa'$. Finally if $I \in \rclass{I}$, it is a
well-known fact that every chain map into $I$ is null homotopic.
Therefore $\ldgclass{A}$ is simply the class of all chain
complexes in Ch($\class{G}$). It is also clear that $\lclass{A} =
\class{E}$, so condition (4) from Theorem~\ref{them-model
structure} holds too. The conclusion of Theorem~\ref{them-model
structure} translates to the usual injective model structure where
the cofibrations are the monomorphisms and the fibrations are the
epimorphisms with dg-injective kernels.
\end{proof}

It was in~\cite{enochs-kaplansky-classes} that the term
``Kaplansky class'' first appeared. The reasoning is that the
class of projective modules (over a ring $R$) form a Kaplansky
class due to Theorem~1 of~\cite{kaplansky}. The details in
interpreting Kaplansky's Theorem in terms of our
Definition~\ref{def-kaplansky class} are in the proof of the next
corollary.

\begin{corollary}\label{cor-projective model structure}
Let $\cat{G} = R$-\textnormal{Mod} where $R$ is a commutative ring
with 1 and let $\class{P}$ be the class of projective modules.
Then $\class{P}$ is a Kaplansky class and it induces the usual
projective model structure on $\textnormal{Ch}(\cat{G})$.
Furthermore this model structure is cofibrantly generated and
monoidal with respect to the usual tensor product of chain
complexes.
\end{corollary}

\begin{proof}
First we note that $\cat{G}$ is locally finitely presentable. For
a proof of this fact, see the footnote to Theorem~4.34
of~\cite{lam}. So $\cat{G}$ is locally $\kappa$-presentable for
every regular cardinal $\kappa$ by Fact~\ref{fact-c.5}.

We now show $\class{P}$ is a Kaplansky class. Since projective
modules are direct summands of free modules, it follows
immediately from Theorem~1 of~\cite{kaplansky} that each
projective module is a direct sum of countably generated modules.
We let $\kappa$ be a regular cardinal with $\kappa
> \text{max} \{ \, |R| , \omega \, \}$ and argue that $\class{P}$ is a $\kappa$-Kaplansky class.
So let $X \subseteq P \neq 0$ where $X$ is $\kappa$-generated and
$P$ is projective. Using Lemmas~\ref{kappa-generated objects
characterized} and~\ref{kappa-generated module is
kappa-presentable module} we are done if we can find a direct
summand $S$ of $P$ which contains $X$ and for which $|S| <
\kappa$. We start by writing $P = \bigoplus_{i \in I} S_i$ where
each $S_i$ is a countably generated (hence $\kappa$-generated)
summand. Let $J = \{ \, j \in I : S_j \cap X \neq 0 \, \}$.
Clearly $|J| < \kappa$ since $|X| < \kappa$. Therefore $S =
\bigoplus_{j \in J} S_j$ is the desired summand since $|S| <
\kappa$.

With the choice of $\kappa$ as above we will now check parts (1) -
(4) of Theorem~\ref{them-model structure}. We just verified
property (1) and property (2) is clearly true since $R \in
\class{P}$ is $\kappa$-generated. Property (3) is true by
Lemma~6.2 of~\cite{hovey} since $\class{P}$ is the left side of a
cotorsion pair. For (4) we use Corollary~3.13 of~\cite{gillespie}.
Basic facts of projective modules show that $\class{P}$ is a
resolving class, that is, $\class{P}$ is closed under taking
kernels of epimorphisms. Referring to Definition~3.11
of~\cite{gillespie} this says that the cotorsion pair
$(\class{P},\class{A})$ (where $\class{A}$ is the class of all
modules) is hereditary. Corollary~3.13 of~\cite{gillespie} now
tells us that $\lclass{P} = \ldgclass{P} \cap \class{E}$. Thus we
have a cofibrantly generated model structure on $\ch$ where the
cofibrations are the monomorphisms with dg-projective cokernels
and the fibrations are the epimorphisms.

Now using Theorem~\ref{them-monoidal condition} it is simple to
check that the model structure is monoidal. It is a standard fact
that projective modules are flat. The tensor product of two
projectives is again projective since $\RHom(P_1 \tensor_R P_2, -)
\cong \RHom(P_1 , \RHom(P_2 -)$ is an exact functor if
$\RHom(P_1,-)$ and $\RHom(P_2,-)$ are exact. Finally, the unit for
the tensor product is $R$ which is projective.
\end{proof}

Another monoidal model structure on $\ch$ which is induced by a
Kaplansky class is the flat model structure which first appeared
in~\cite{gillespie}.

\begin{corollary}\label{cor-flat model structure}
Let $\cat{G} = R$-\textnormal{Mod} where $R$ is a commutative ring
with 1 and let $\class{F}$ be the class of flat modules. Then
$\class{F}$ is a Kaplansky class and it induces a model structure
on \textnormal{Ch}($\cat{G}$) called the flat model structure. The
cofibrations are the monomorphisms with dg-flat cokernels and the
fibrations are the epimorphisms with dg-cotorsion kernels.
Furthermore this model structure is cofibrantly generated and
monoidal with respect to the usual tensor product of chain
complexes.
\end{corollary}

\begin{proof}
Again we note that $\cat{G}$ is locally finitely presentable. So
it is locally $\kappa$-presentable for every regular cardinal
$\kappa$. We let $\kappa$ be a regular cardinal with $\kappa
> |R|$ and argue that $\class{F}$ is a $\kappa$-Kaplansky class.

So let $X \subseteq F \neq 0$ where $X$ is $\kappa$-generated and
$F$ is flat. By Lemma~\ref{kappa-generated objects characterized},
$|X| < \kappa$. For each $x \in X$ we use Lemma~2
of~\cite{enochs-flat-cover-theorem} to find a flat submodule $F_x
\subseteq F$ with $x \in F_x$, $|F_x| < \kappa$ and $F/F_x \in
\class{F}$. The direct union $\bigcup_{x \in X} F_x$ contains $X$
and is flat. Furthermore $$F/(\bigcup_{x \in X} F_x) \cong \lim_{x
\in X} F/F_x$$ is flat. Finally, since $|\bigcup_{x \in X} F_x| <
\kappa$, Lemmas~\ref{kappa-generated objects characterized}
and~\ref{kappa-generated module is kappa-presentable module} allow
us to conclude that $\class{F}$ is a $\kappa$-Kaplansky class.

This proves property (1) of Theorem~\ref{them-model structure} and
property (2) is clearly true since $R$ is a flat $R$-module.
Proposition~XVI.3.1 of~\cite{lang} tells us that $\class{F}$ is
closed under retracts (direct summands). It is also a standard
fact about flat modules that $\class{F}$ is closed under
extensions and direct limits. It follows that $\class{F}$ is
closed under transfinite extensions. So it is left to check
property (4) of Theorem~\ref{them-model structure}. By
Proposition~XVI.3.4 of~\cite{lang} $\class{F}$ is closed under
taking kernels of epimorphisms. Therefore we can argue as in the
proof of Corollary~\ref{cor-projective model structure} that
$\lclass{F} = \ldgclass{F} \cap \class{E}$.

Thus we have a cofibrantly generated model structure on $\ch$
where the cofibrations are the monomorphisms with cokernels in
$\ldgclass{F}$ (dg-flat complexes) and the fibrations are the
epimorphisms with kernels in $\rdgclass{C}$ (dg-cotorsion
complexes).

Again it is easy to see that the model structure is monoidal using
Theorem~\ref{them-monoidal condition}. The tensor product of two
flat modules is clearly flat and the unit $R$ is flat.
\end{proof}

As described in~\cite{gillespie-sheaves} the (monoidal) flat model
structure of Corollary~\ref{cor-flat model structure} generalizes
to the category Ch($\class{O}_X$-Mod) of complexes of
$\class{O}_X$-modules where $\class{O}_X$ is a ringed space. It
too is obtained from a Kaplansky class as we will see next in
Corollary~\ref{cor-flat model structure on sheaves}. The first
lemma below concerns cotorsion modules and skyscraper sheaves. We
now recall the concept of a skyscraper sheaf.

Let $\class{O}_X$ be a ringed space, $p \in X$ be a point, and $M$
be an $\class{O}_p$-module. The skyscraper sheaf $S_p (M)$ is the
sheaf on $X$ defined by $U \mapsto M$ if $p \in U$ and $U \mapsto
0$ if $p \notin U$. It is in fact an $\class{O}_X$-module by
viewing $M$ as an $\class{O}(U)$-module via the ring homomorphism
$\class{O}(U) \ar \class{O}_p$. One can check that $[S_p (M)]_q =
M$ for each $q \in \overline{\{\, p \, \}}$ and $[S_p (M)]_q = 0$
for each $q \notin \overline{\{\, p \, \}}$. It is also standard
that $S_p (-)$ is an exact functor from the category
$\class{O}_p$-Mod to the category $\class{O}_X$-Mod and is right
adjoint to the (also exact) ``stalk functor'' which sends a
$\class{O}_X$-module $F$ to the $\class{O}_p$-module $F_p$.
Therefore $S_p(-)$ preserves injective objects by
Proposition~2.3.10 of~\cite{weibel}.

\begin{lemma}\label{lemma-skyscraper functor preserves cotorsion}
Let $\class{O}_X$ be a ringed space and $p \in X$ be a point. The
skyscraper functor $S_p(-)$ preserves cotorsion objects.
\end{lemma}

\begin{proof}
First notice that for any $\class{O}_X$-module $F$ and any
$\class{O}_p$-module $C$ we have an isomorphism $\Ext^n(F_p,C)
\cong \Ext^n(F,S_p(C))$. (This can be proved using the adjoint
relationship discussed above along with the fact that $S_p(-)$ is
exact and preserves injective objects, and therefore preserves
injective resolutions.)

Now suppose $C$ is a cotorsion $\class{O}_p$-module and $F$ is a
flat $\class{O}_X$-module. We wish to show that $S_p(C)$ is a
cotorsion $\class{O}_X$-module. But $F_p$ is a flat
$\class{O}_p$-module, so $\Ext^1(F,S_p(C)) \cong \Ext^1(F_p,C) =
0$. Therefore $S_p(C)$ is cotorsion.
\end{proof}

In order to verify the first hypothesis of Theorem~\ref{them-model
structure}, we would like to find a large enough cardinal $\kappa$
so that the $\kappa$-generated $\class{O}_X$-modules coincide with
the $\kappa$-presentable $\class{O}_X$-modules and so that these
notions may be used in place of cardinality.
Lemma~\ref{lemma-cardinality of sheaf and large kappas} will allow
us to do this. Our next lemma below however states that the
category of $\class{O}_X$-modules has a set of flat generators.
This is a standard fact but we document it now since it will be
used in the proof of Lemma~\ref{lemma-cardinality of sheaf and
large kappas} and is needed for Corollary~\ref{cor-flat model
structure on sheaves}.

Let $\class{O}_X$ be a ringed space. For each open $U \subseteq
X$, extend $\class{O}_{|U}$ by 0 outside of $U$ to get a presheaf,
denoted $\class{O}_U$. Now sheafify to get an
$\class{O}_X$-module, which we will denote $j!(\class{O}_U)$.

\begin{lemma}\label{lemma-set of flat generators}
Let $\class{O}_X$ be a ringed space. Then $\Hom(j!(\class{O}_U) ,
G) \cong G(U)$ for any $\class{O}_X$-module $G$ and each open set
$U \subseteq X$. Furthermore, $\{ \, j!(\class{O}_U) : U \subseteq
X \, \}$ is a set of flat generators for
$\class{O}_X$-\textnormal{Mod}.
\end{lemma}

\begin{proof}
One can prove without much difficulty that $\Hom(\class{O}_U , G)
\cong G(U)$. So by the universal property of sheafification we get
$\Hom(j!(\class{O}_U) , G) \cong G(U)$. It follows at once that
the set forms a generating set since the modules $j!(\class{O}_U)$
``pick out points''. Also, each $j!(\class{O}_U)$ is flat since
$[j!(\class{O}_U)]_p \cong (\class{O}_U)_p$, which equals
$\class{O}_p$ if $p \in U$ and  0 if $p \in X \backslash U$.
\end{proof}

\begin{definition}
We define the cardinality of a presheaf (or sheaf), $F$, to be
$|F| = |\coprod_{U \subseteq X} F(U)|$ where $U \subseteq X$
ranges over all the open sets in $X$.
\end{definition}

In Lemmas~\ref{lemma-sheafification has same cardinality for large
kappas} and~\ref{lemma-cardinality of sheaf and large kappas} we
let $\class{O}_X$ be our ringed space and let $\class{U}$
represent the set of all open sets $U \subseteq X$. We note now
that if $\beta$ is an infinite cardinal in which $\beta >
|\class{O}_X|$, then automatically $\beta > |\class{U}|$, since
$\class{O}(U)$ is nonempty for each $U \subseteq X$.

\begin{lemma}\label{lemma-sheafification has same cardinality for
large kappas}
Let $\beta$ be an infinite cardinal such that $\beta
> \textnormal{max} \{ \, |X| , |\class{O}_X| \, \}$.
Now let $\kappa = 2^{\beta}$. If $S$ is a presheaf of
$\class{O}_X$-modules and $|S| < \kappa$, then $|S^+| < \kappa$
where $S^+$ is the sheafification.
\end{lemma}

\begin{proof}
Clearly $|S_p| < \kappa$ for each point $p \in X$. So by the
sheafification construction (see the proof of
Proposition-Definition~II.1.2 of~\cite{hartshorne}) one can easily
see that for each open $U \subseteq X$ we have $|S^+(U)| <
\kappa^{\beta}$. However, $\kappa^{\beta} = (2^{\beta})^{\beta} =
2^{(\beta^2)} = 2^{\beta} = \kappa$, so $|S^+(U)| < \kappa$.
Therefore $|S^+| = |\coprod_{U \subseteq X} S^+(U)| < \kappa$.
\end{proof}

Next suppose we have a presheaf of $\class{O}_X$-modules, $S$, and
suppose we have a subset $$W \subseteq \coprod_{U \subseteq X}
S(U)$$ where again $U \subseteq X$ ranges over all the open sets
in $X$. We set $W_U = W \cap S(U)$. The set $W$ generates a
presheaf $S'$ and it is given by
$$S'(U) = \sum_{V \supseteq U} \, \sum_{w \in W_V} r_{V,U} ([\class{O}(V)]w).$$
That is, $S'(U)$ consists of all finite sums of the form
$$r_{V_1,U}(\rho_1 w_1) + r_{V_2,U}(\rho_2 w_2) + \cdots +
r_{V_n,U}(\rho_n w_n)$$ where $$V_i \supseteq U, \, w_i \in
W_{V_i}, \, \rho_i \in \class{O}(V_i)$$ and $r_{V_i,U} \mathcolon
V_i \ar U$ are the restriction maps. It is straightforward to see
that $S'$ is a presheaf and upon reflection it is clearly the
smallest subpresheaf containing $W$. Furthermore, if $\kappa
> |\class{O}_X|$ is an infinite cardinal, and $\kappa > |W|$, then
we will have $\kappa > |S'|$. If $S$ was a sheaf of
$\class{O}_X$-modules to start with, then $(S')^+$ is the
$\class{O}_X$-submodule generated by $W$ which we will denote
$S_W$. In this case note that if $\kappa = 2^{\beta}$ where $\beta
> \textnormal{max} \{ \, |X| , |\class{O}_X| \,
\}$, then by Lemma~\ref{lemma-sheafification has same cardinality
for large kappas},  $\kappa > |S_W|$ whenever $\kappa > |W|$.

\begin{lemma}\label{lemma-cardinality of sheaf and large kappas}
Let $\beta$ be an infinite cardinal such that $\beta
> \textnormal{max} \{ \, |X| , |\class{O}_X| \, \}$.
Now let $\kappa = 2^{\beta}$. Also assume that $\kappa$ is large
enough that each $j!(\class{O}_U)$ is $\kappa$-generated. Then the
following are equivalent for an $\class{O}_X$-module $S$.

(1) $|S| < \kappa$.

(2) $S$ is $\kappa$-generated.

(3) $S$ is $\kappa$-presentable.

\end{lemma}

\begin{proof}
(1) $\Rightarrow (2)$. We use
Fact~\ref{kappa-generated-characterization}. Suppose $|S| <
\kappa$ and say $S = \sum_{i \in I} S_i$ is a $\kappa$-directed
union of $\class{O}_X$-submodules. Note that we will be done if we
can show that for each open $U \subseteq X$ and each $x \in S(U)$,
there exists an $i \in I$ such that $x \in S_i(U)$. Indeed if this
were true, then by the large choice of $\kappa$ and the fact that
the union is $\kappa$-filtered, we would take the union of all
such $S_i(U)$ to display $S = S_i$ for some particular $i \in I$.
(Note also that the assertion we wish to prove is not obvious
since we must sheafify when taking the direct union.)

Now let $x \in S(U)$. From our choice of $\kappa$ we know that
each $j!(\class{O}_U)$ is $\kappa$-generated and so the canonical
map $\colim_{i \in I} \Hom(j!(\class{O}_U) , S_i) \ar
\Hom(j!(\class{O}_U) , S)$ is an isomorphism. Now using
Lemma~\ref{lemma-set of flat generators} this translates to an
isomorphism $\colim_{i \in I} S_i(U) \cong S(U)$. Through this
isomorphism we see that $x \in S_i(U)$ for some $i \in I$.

$(2) \Rightarrow (1)$. Let $\class{W}$ be the collection of
\emph{all} subsets $W \subseteq \coprod_{U \subseteq X} S(U)$
(where $U$ ranges over all open subsets of $X$) which satisfy $|W|
< \kappa$. For each $W \in \class{W}$, let $S_W$ represent the
$\class{O}_X$-submodule generated by $W$. Then $|S_W| < \kappa$ by
Lemma~\ref{lemma-sheafification has same cardinality for large
kappas}. Note that $(\class{W}, \subseteq)$ is $\kappa$-filtered
and in fact $S$ is the $\kappa$-filtered union $\sum_{W \in
\class{W}} S_W$. By Fact~\ref{kappa-generated-characterization},
$S = S_W$ for some $W \in \class{W}$. So $|S| < \kappa$.

$(3) \Rightarrow (2)$ is automatic. We now prove $(2) \Rightarrow
(3)$, using that (1) iff (2). First we point out that the category
of $\class{O}_X$-modules is locally $\kappa$-generated. Indeed
each $S$ can be expressed as the $\kappa$-filtered union $S =
\sum_{W \in \class{W}} S_W$ where each $S_W$ is $\kappa$-generated
as in the last paragraph. Therefore, we may use the
characterization of $\kappa$-presentable objects in
Fact~\ref{kappa-prentable-characterization}. Suppose $S$ is
$\kappa$-generated and $T \ar S$ is an epimorphism with $T$ a
$\kappa$-generated $\class{O}_X$-module. Then $|T| < \kappa$. So
of course $|\ker{(T \ar S)}| < \kappa$, which means $\ker{(T \ar
S)}$ is $\kappa$-generated. This proves $S$ is
$\kappa$-presentable.
\end{proof}

Finally we prove that the flat model structure on
Ch($\class{O}_X$-Mod) comes from the Kaplansky class of flat
$\class{O}_X$-modules. The proof of Corollary~\ref{cor-flat model
structure on sheaves} will again rely on a result of Enochs,
et.al. which can be found in~\cite{enochs-oyonarte}.

If $S \subseteq F$ is a subpresheaf (or subsheaf) we call $S$
\emph{presheaf pure} if $S(U)$ is a pure $\class{O}(U)$-submodule
of $F(U)$ for each open $U$. We say $S \subseteq F$ is
\emph{stalkwise pure} if $S_p$ is a pure $\class{O}_p$-submodule
of $F_p$ for each $p \in X$.

\begin{corollary}\label{cor-flat model structure on sheaves}
Let $\cat{G} = \class{O}_X$-\textnormal{Mod} where $\class{O}_X$
is a sheaf of rings on a topological space $X$ and let $\class{F}$
be the class of flat $\class{O}_X$-modules. Then $\class{F}$ is a
Kaplansky class and it induces a model structure on
\textnormal{Ch}($\cat{G}$) we call the flat model structure. The
cofibrations are the monomorphisms with dg-flat cokernels and the
fibrations are the epimorphisms with dg-cotorsion kernels.
Furthermore this model structure is cofibrantly generated and
monoidal with respect to the usual tensor product of chain
complexes.
\end{corollary}

\begin{proof}
Let $\mathcal{U}$ represent the set of all open sets $U \subseteq
X$. Let $\beta$ be an infinite cardinal such that $\beta
> \textnormal{max} \{ \, |X| , |\mathcal{U}| , |\class{O}_X| \, \}$.
Now let $\kappa = 2^{\beta}$.  Using Fact~\ref{fact-a} we can also
assume that $\kappa$ is large enough that each $j!(\class{O}_U)$
is $\kappa$-generated. We let $\class{F}$ be the class of flat
sheaves and claim that $\class{F}$ is a $\kappa$-Kaplansky class.

So let $S' \subseteq F \neq 0$ where $S'$ is $\kappa$-generated
and $F$ is flat. By Lemma~\ref{lemma-cardinality of sheaf and
large kappas}, $|S'| < \kappa$. With minor adjustments to the
proof of Proposition~2.4 of~\cite{enochs-oyonarte} we can find an
$\class{O}_X$-submodule $S \subseteq F$ which is presheaf pure and
such that $S' \subseteq S$ and $|S| < \kappa$. It follows that $S
\subseteq S^+$ and $S^+ \subseteq F$ is stalkwise pure. (This is
not hard, but one could also see a proof in Lemma~4.7
of~\cite{gillespie-sheaves}.) It follows immediately that $S^+$
and $F/S^+$ are flat. Furthermore, by
Lemma~\ref{lemma-sheafification has same cardinality for large
kappas} we have $|S^+| < \kappa$ and so by
Lemma~\ref{lemma-cardinality of sheaf and large kappas} we have
that $S^+$ is $\kappa$-presentable. This completes the proof that
$\class{F}$ is a $\kappa$-Kaplansky class.

Note that $\class{O}_X$-Mod is locally $\kappa$-presentable by the
same argument given in the last paragraph of the proof of
Lemma~6.8. So we are done verifying properties (1) and (2) of
Theorem~\ref{them-model structure}.

It is a standard fact that $\class{F}$ is closed under extensions
and direct limits. It follows that $\class{F}$ is closed under
transfinite extensions. It is also standard that $\class{F}$ is
closed under retracts. So we focus on proving property (4). We let
$Y$ be an exact dg-flat complex of sheaves. We want to show that
it is a flat complex. I.e. that is has flat cycles in each degree.
Note that we will be done if we can show that the``stalk complex''
$Y_p$ is a flat complex of $\class{O}_p$-modules for each $p \in
X$. We know from the proof of Corollary~\ref{cor-flat model
structure} that an exact dg-flat complex of $\class{O}_p$-modules
is a flat complex. So we will simply show that $Y_p$ is a dg-flat
complex. So let $Y_p \ar C$ be a morphism where $C$ is a cotorsion
complex of $\class{O}_p$-modules. (This means that $C$ is exact
and each cycle is cotorsion.) We want to show this map is null
homotopic. But using the skyscraper functor and its adjoint
property gives us a morphism $Y \ar S_p(C)$ of complexes. Since
$S_p(-)$ is exact and preserves cotorsion modules by
lemma~\ref{lemma-skyscraper functor preserves cotorsion} we see
that $S_p(C)$ is a cotorsion complex and by definition of dg-flat
$Y \ar S_p(C)$ is null homotopic. From this we can deduce that
$Y_p \ar C$ is also null homotopic.

Thus we have a cofibrantly generated model structure on
Ch($\cat{G})$ where the cofibrations are the monomorphisms with
cokernels in $\ldgclass{F}$ (dg-flat complexes) and the fibrations
are the epimorphisms with kernels in $\rdgclass{C}$ (dg-cotorsion
complexes).

We check now that the model structure is monoidal using
Theorem~\ref{them-monoidal condition}. Again this is easy. The
tensor product of two flat sheaves is also flat. The unit is the
structure sheaf $\class{O}_X$ which is flat too.
\end{proof}

\appendix

\section{locally presentable
categories}\label{appendix-locally-presentable-cats}

Tibor Beke showed (Proposition~3.10 of~\cite{beke}) that
Grothendieck categories are locally presentable. Here we gather
some properties of locally presentable categories from the
literature which will be needed for this paper. Virtually
everything here can be found scattered
throughout~\cite{adamek-rosicky} and~\cite{stenstrom}.

Throughout this section, we will assume all cardinals are regular.
These are infinite cardinals which are not the sum of a smaller
number of smaller cardinals. For example, $\aleph_{\omega} =
\sum_{n < \omega} \aleph_n$ is NOT a regular cardinal. All
infinite successor cardinals are regular though.

\begin{definition} Let $\cat{C}$ be a category and let $\kappa$ be a regular cardinal.

A \emph{$\kappa$-filtered category} is a category $\cat{K}$, for
which every subcategory with less than $\kappa$ morphisms has a
cocone. A \emph{$\kappa$-filtered diagram} is simply a functor $F
\mathcolon \cat{K} \ar \cat{C}$ in which $\cat{K}$ is a small
$\kappa$-filtered category. By a \emph{$\kappa$-filtered colimit},
we mean the colimit of a $\kappa$-filtered diagram.

An object $X$ in $\cat{C}$ is called \emph{$\kappa$-presentable}
if $\Hom_{\cat{C}}(X,-)$ preserves $\kappa$-filtered colimits. It
is called \emph{presentable} if it is $\kappa$-presentable for
some regular cardinal $\kappa$.

An object $X$ in $\cat{C}$ is called \emph{$\kappa$-generated} if
$\Hom_{\cat{C}}(X,-)$ preserves $\kappa$-filtered colimits of
monomorphisms. I.e., preserves the colimits of diagrams $F
\mathcolon \cat{K} \ar \cat{C}$ for which $F(d) \mathcolon F(c)
\ar F(c')$ is a monomorphism for each $d \mathcolon c \ar c'$ in
$\cat{K}$. The object is called \emph{generated} if it is
$\kappa$-generated for some regular cardinal $\kappa$.
\end{definition}

When $\cat{C}$ is the category of $R$-modules, the definitions
above for $\kappa$-presentable and $\kappa$-generated agree with
usual conventions. In fact, see below
Fact~\ref{kappa-prentable-characterization} and the example after
Fact~\ref{kappa-generated-characterization}.

As a special case we have the notion of a \emph{$\kappa$-directed
poset} which is really just a poset $(P,\leq)$ which is
$\kappa$-filtered when thought of as a category. Then we have
\emph{$\kappa$-directed colimits} which are colimits of a functor
$F \mathcolon (P,\leq) \ar \cat{K}$, where $(P,\leq)$ is a
$\kappa$-directed poset. In a Grothendieck category $\cat{G}$, we
often use the term \emph{$\kappa$-directed union}. This is just a
$\kappa$-directed colimit of a $\kappa$-directed set of
subobjects, $\{ S_i \}_{i \in I}$, of a given object $A$. It is
denoted $\sum_{i \in I} S_i$ and it is a subobject of $A$ since
direct limits are exact in a Grothendieck category.

\begin{fact}\label{fact-a}
Let  $\kappa' \geq \kappa$ be regular cardinals. It is easy to see
that any $\kappa'$-filtered category is also a $\kappa$-filtered
category. Thus any $\kappa'$-filtered diagram is a
$\kappa$-filtered diagram. Therefore, any $\kappa$-presentable
(resp. generated) object is easily seen to be
$\kappa'$-presentable (resp. generated).
\end{fact}

\begin{fact}\label{fact-b}
Clearly, any $\kappa$-presentable object is $\kappa$-generated.
\end{fact}

As explained in~\cite{adamek-rosicky}, $\kappa$-presentable and
$\kappa$-generated objects can be defined in terms of
$\kappa$-directed colimits instead of $\kappa$-filtered colimits.
We also have the following facts which hold by generalizing the
proofs in Chapter~V.3 of~\cite{stenstrom}. The proofs there are
given for $\kappa = \omega$ (finitely generated objects), but
everything carries over for an arbitrary $\kappa$.

\begin{fact}\label{fact-b.5}
In a Grothendieck category, the image of a $\kappa$-generated
object is again $\kappa$-generated. Also an extension of
$\kappa$-generated objects is $\kappa$-generated.
\end{fact}

\begin{fact}\label{kappa-generated-characterization}
Let $\cat{G}$ be a Grothendieck category and $\kappa$ a regular
cardinal. Then the following are equivalent for an object $X \in
\cat{G}$:

1) $X$ is $\kappa$-generated.

2) Whenever $X = \sum_{i \in I} X_i$ is a $\kappa$-directed union
of subobjects, we have $X = X_i$ for some $i \in I$.

3) Whenever $X = \sum_{i \in I} X_i$ is \emph{any} union of
subobjects, we have $X = \sum_{i \in J} X_i$ for some $J \subseteq
I$ with $|J| < \kappa$.
\end{fact}

\begin{proof}
The first two are equivalent by generalizing the proofs in
Chapter~V.3 of~\cite{stenstrom}.  The equivalence of the second
two is not hard to prove but you need the fact that a colimit of
less than $\kappa$-many $\kappa$-generated objects is again
$\kappa$-generated.
\end{proof}

\begin{example}\label{kappa-generated for modules}
Let $R$ be a ring and $M$ an $R$-module. $M$ is $\kappa$-generated
iff there exists a set $S \subseteq M$ with $|S| < \kappa$ such
that $M = \sum_{x \in S} Rx$.

\begin{proof}
$(\Rightarrow)$. Let $\class{S}$ be the collection of \emph{all}
subsets $S \subseteq M$ such that $|S| < \kappa$. Let $M_S$ be the
submodule $\sum_{x \in S} Rx$. Note that $(\class{S}, \subseteq)$
is $\kappa$-filtered and in fact $M$ is the $\kappa$-filtered
union $\sum_{S \in \class{S}} M_S$. Using
Fact~\ref{kappa-generated-characterization}, $M = M_S$ for some $S
\in \class{S}$.

($\Leftarrow)$. Again we will use
Fact~\ref{kappa-generated-characterization}. Let $M = \sum_{i \in
I} M_i$ be a $\kappa$-filtered union. By assumption, there is an
$S \subseteq M$ such that $M = M_S \ (|S| < \kappa)$. For each $x
\in S$, there corresponds some $i \in I$ such that $x \in M_i$.
But $\sum_{i \in I} M_i$ is $\kappa$-filtered and $|S| < \kappa$,
so there exists $i_0 \in I$ such that $S \subseteq M_{i_0}$. So $M
= M_{i_0}$.
\end{proof}

\end{example}

\begin{definition}
Any cocomplete category $\cat{C}$ is called \emph{locally
$\kappa$-presentable} if each object is a $\kappa$-filtered
colimit of $\kappa$-presentable objects and the class of
$\kappa$-presentable objects is \emph{essentially small}, meaning
there are only a sets worth of $\kappa$-presentable objects up to
isomorphism. A category $\cat{C}$ is simply called \emph{locally
presentable} if it is locally $\kappa$-presentable for some
regular cardinal $\kappa$.

Analogously, we define a \emph{locally $\kappa$-generated}
category as a cocomplete category in which each object is a
$\kappa$-filtered union of $\kappa$-generated subobjects and the
class of $\kappa$-generated objects is essentially small. A
category $\cat{C}$ is simply called \emph{locally generated} if it
is locally $\kappa$-generated for some regular cardinal $\kappa$.
\end{definition}

For example, if $R$ is a ring, then every $R$-module is an
$\omega$-filtered union (usually just called filtered union) of
its $\omega$-generated submodules (finitely generated submodules).
So the category of $R$-modules is locally $\omega$-generated
(locally finitely generated).

In general, as pointed out in~\cite{stenstrom} (pp. 122), an
arbitrary Grothendieck category may not even have  any nonzero
finitely generated objects. Such a category is clearly not locally
finitely generated. However, the ``Local Generation
Theorem''(~\cite{adamek-rosicky}, pp.54) says that a cocomplete
category $\cat{C}$ is locally presentable if and only if it is
locally generated. But be careful; this is NOT the same as saying
$\cat{C}$ is locally $\kappa$-presentable if and only if it is
locally $\kappa$-generated. For an example,
see~\cite{adamek-rosicky}, Example~1.71, pp.56. Nevertheless, one
direction does hold. For a proof of the following fact, see the
first paragraph of the proof to the ``Local Generation
Theorem''(~\cite{adamek-rosicky}, pp.54).

\begin{fact}\label{fact-c}
If $\cat{C}$ is locally $\kappa$-presentable, then $\cat{C}$ is
locally $\kappa$-generated.
\end{fact}

\begin{fact}\label{fact-c.5}
Let $\kappa' \geq \kappa$ be regular cardinals. Then $\cat{C}$
locally $\kappa$-presentable (resp. generated) implies $\cat{C}$
locally $\kappa'$-presentable (resp. generated).
\end{fact}

\begin{fact}\label{fact-d}
If $\cat{C}$ is locally presentable, then for any regular cardinal
$\lambda$, the class of $\lambda$-presentable objects is
essentially small. Similarly, the class of $\lambda$-generated
objects is essentially small(\cite{adamek-rosicky},
Corollary~1.69). Any fixed set of representatives of the
$\lambda$-presentable objects is denoted $\textbf{Pres}_{\lambda}
\ \cat{C}$, and a set of representatives of the
$\lambda$-generated objects is denoted by $\mathbf{Gen}_{\lambda}
\ \cat{C}$.
\end{fact}

\begin{fact} Every object of a locally presentable category is
presentable (and hence generated).
\end{fact}

\begin{fact}\label{fact-wellpowered} Let $X$ be an object in a locally presentable
category. Then the class of subobjects of $X$ is in fact a set and
the class of quotient objects of $X$ is also a set.
In~\cite{adamek-rosicky} this is expressed by saying locally
presentable categories are \emph{wellpowered} and
\emph{cowellpowered}, respectively.
\end{fact}

The following gives a useful characterization of
$\kappa$-presentable objects. Again, the proof is obtained by
generalizing the proofs in Chapter~V.3 of~\cite{stenstrom}.

\begin{fact}\label{kappa-prentable-characterization}
Let $\cat{G}$ be a locally $\kappa$-generated Grothendieck
category. An object $X$ is $\kappa$-presentable if and only if $X$
is $\kappa$-generated and every epimorphism $C \ar X$  with $C$
$\kappa$-generated has a $\kappa$-generated kernel.
\end{fact}


\section{$\kappa$-generated objects for very large
$\kappa$}\label{appendix-kappa-generated and kappa-presentable
objects coincide}

Here we prove that in a Grothendieck category, there exists a
regular cardinal $\kappa$ for which the class of
$\kappa$-generated objects coincides with the class of
$\kappa$-presentable objects. As a result, this class of objects
(which has a small skeleton) satisfies the 2 out of 3 property for
short exact sequences.

\begin{lemma}\label{kappa-generated objects characterized}
Let $R$ be a ring and $M$ an $R$-module. Let $\kappa$ be a regular
cardinal such that $\kappa > |R|$. Then $M$ is $\kappa$-generated
iff $|M| < \kappa$.
\end{lemma}

\begin{proof}
We use the characterization of $\kappa$-generated provided by
Example~\ref{kappa-generated for modules}.

First, if $|M| < \kappa$, then $M$ is clearly $\kappa$-generated
since we can take the generating set to be $M$ itself. Conversely,
suppose $M$ is $\kappa$-generated so that there exists a set $S
\subseteq M$ for which $M = \sum_{x \in S} Rx$ and $|S| < \kappa$.
Then $|M| = |\sum_{x \in S} Rx| \leq |\bigoplus_{x \in S} Rx| <
\kappa$.

\end{proof}

\begin{lemma}\label{kappa-generated module is kappa-presentable
module} Let $R$ be a ring and $M$ an $R$-module. Let $\kappa$ be a
regular cardinal such that $\kappa > |R|$. Then $M$ is
$\kappa$-generated iff $M$ is $\kappa$-presentable.
\end{lemma}

\begin{proof}
Of course $\kappa$-presentable objects are always
$\kappa$-generated, so it remains to show that under the given
hypothesis, $\kappa$-generated objects are $\kappa$-presentable.
But this is easy using Lemma~\ref{kappa-generated objects
characterized} and the fact that $M$ is $\kappa$-presentable iff
$M$ is $\kappa$-generated and every epimorphism $N \ar M$ with $N$
$\kappa$-generated, has a $\kappa$-generated kernel. Indeed
suppose $M$ is $\kappa$-generated and $N \ar M$ is an epimorphism
with $N$ $\kappa$-generated. Then $|N| < \kappa$. So of course
$|\ker{(N \ar M)}| < \kappa$.
\end{proof}

Next we can prove that the notion of $\kappa$-generated coincides
with $\kappa$-presentable in any Grothendieck category if we take
$\kappa$ to be large enough. In the proof we use the fact that an
equivalence of categories preserves $\kappa$-generated and
$\kappa$-presentable objects. Indeed if $F \mathcolon \cat{A} \ar
\cat{B}$ is an equivalence, then there exists $G \mathcolon
\cat{B} \ar \cat{A}$ such that $GF \cong 1_{\cat{A}}$ and $FG
\cong 1_{\cat{B}}$. Furthermore, $F$ is then both a left and right
adjoint of $G$ (and so each preserves colimits). So if $X \in
\cat{A}$ is $\kappa$-presentable then the functor
$\cat{A}(X,G(-))$ preserves $\kappa$-filtered colimits. But
$\cat{A}(X,G(-)) \cong \cat{B}(F(X),-)$, so $F(X)$ must be
$\kappa$-presentable. On the other hand, if $F(X)$ is
$\kappa$-presentable, then $\cat{B}(F(X),F(-))$ preserves
$\kappa$-filtered colimits. Therefore, $\cat{B}(F(X),F(-)) \cong
\cat{A}(GF(X),-)$ preserves $\kappa$-filtered colimits and $GF(X)
\cong X$ must be $\kappa$-presentable. The same argument (but with
$\kappa$-filtered colimits of monomorphisms) shows that an
equivalence preserves $\kappa$-generated objects.

\begin{proposition}
Let $\cat{G}$ be a Grothendieck category. Then there exists a
regular cardinal $\kappa$ for which the $\kappa$-generated objects
coincide with the $\kappa$-presentable objects.
\end{proposition}

\begin{proof}
By the Gabriel-Popescu Theorem (\cite{stenstrom}), $\cat{G}$ is
equivalent to a subcategory of $R$-Mod for some ring $R$. (The
Gabriel-Popescu Theorem says more than this, but we don't need the
whole statement.) Since an equivalence of categories preserves
$\kappa$-generated and $\kappa$-presentable objects, the result
follows from Lemma~\ref{kappa-generated module is
kappa-presentable module} by choosing $\kappa$ to be a regular
cardinal $\kappa
> |R|$.
\end{proof}


\end{document}